\documentclass[noams,lp]{compositio}
\usepackage{amsbsy,amssymb,amscd,amsfonts,latexsym,amstext,delarray,
 amsmath,stackrel,lineno,tabularx,url,MnSymbol,tikz,cite,yfonts}

\input xypic
\usepackage[pdfborder={0 0 0}]{hyperref}

\hypersetup{
    colorlinks = true,
    linkcolor = blue,
    anchorcolor = red,
    citecolor = green,
    filecolor = red,
    pagecolor = red,
    urlcolor = magenta
}
\newcommand*\patchAmsMathEnvironmentForLineno[1]{%
  \expandafter\let\csname old#1\expandafter\endcsname\csname #1\endcsname
  \expandafter\let\csname oldend#1\expandafter\endcsname\csname end#1\endcsname
  \renewenvironment{#1}%
     {\linenomath\csname old#1\endcsname}%
     {\csname oldend#1\endcsname\endlinenomath}}%
\newcommand*\patchBothAmsMathEnvironmentsForLineno[1]{%
  \patchAmsMathEnvironmentForLineno{#1}%
  \patchAmsMathEnvironmentForLineno{#1*}}%
\AtBeginDocument{%
\patchBothAmsMathEnvironmentsForLineno{equation}%
\patchBothAmsMathEnvironmentsForLineno{align}%
\patchBothAmsMathEnvironmentsForLineno{flalign}%
\patchBothAmsMathEnvironmentsForLineno{alignat}%
\patchBothAmsMathEnvironmentsForLineno{gather}%
\patchBothAmsMathEnvironmentsForLineno{multline}%
}
\newtheorem{thm}{Theorem}[section] % number like 3.1, 3.2, 3.3, etc.
\newtheorem{defn}[thm]{Definition} % numbered with thm
\newtheorem{prop}[thm]{Proposition}

\newtheorem{lem}[thm]{Lemma}

\newtheorem{rem}[thm]{Remark}

\def\Aut{{\rm Aut}}

\def\End{{\rm End}}

\def\GL{{\rm GL}}

\def\Hom{{\rm Hom}}
\def\id{{\rm id}}

\def\mmod{{\rm Mod}}

\def\Spec{{\rm Spec\,}}

\def\Tr{{\rm Tr}}

\def\A{{\mathbb A}}
\def\B{{\mathbb B}}
\def\C{{\mathbb C}}
\def\F{{\mathbb F}}
\def\K{{\mathbb K}}
\def\N{{\mathbb N}}

\def\Q{{\mathbb Q}}
\def\R{{\mathbb R}}
\def\Z{{\mathbb Z}}

\def\urep{\vartheta}

\def\pt{\mathfrak P}

\def\Tr{{\rm Tr}}

\def\cC{{\mathcal C}}

\def\cE{{\mathcal E}}

\def\cH{{\mathcal H}}

\def\cF{{\mathcal F}}

\def\cK{{\mathcal K}}

\def\cM{{\mathcal M}}
\def\cO{{\mathcal O}}

\def\cR{{\mathcal R}}
\def\cS{{\mathcal S}}
\def\cT{{\mathcal T}}

\def\qcy{{\Q^{\rm cyc}}}
\def\ocy{{\cO^{\rm cyc}}}
\def\card{{\mathcal{C}\it{a}\mathcal{C}\ell}}

\def\qqq{\,,\,~\forall}

\def\sub{{\rm Sub}_\geq}
\def\conv{{\rm Conv}_\geq}

\newcommand{\ie}{{\it i.e.\/}\ }
\newcommand{\eg}{{\it e.g.\/}\ }
\newcommand{\cf}{{\it cf.}}
\newcommand{\opcit}{{\it op.cit.\/}\ }

\def\no{\noindent}

\def\spz{{\Spec\Z}}

\def\zmin{\Z_{\rm min}}
\def\zminp{\Z_{\rm min}^+}

\def\id{{\mbox{Id}}}

\def\Hom {{\mbox{Hom}}}

\def\End{{\mbox{End}}}

\def\Int{{\mbox{Int}}}

\def\catmo{{\mathfrak{Mod}}}

\def\ff{{\rm  Fr}}

\def\ffp{\mathfrak{p}}

\def\mc{multiplicatively cancellative }
\def\zmax{{\Z_{\rm max}}}
\def\zhmax{{\Z_{\rm hmax}}}

\def\rmax{\R_+^{\rm max}}

\def\fr{{\rm Fr}}

\def\arith{{(\wnt,\bar \N)}}

\def\Se{\frak{ Sets}}

\def\nt{\N^{\times}}
\def\wnt{{\widehat{\N^{\times}}}}
\def\nto{\N_0^{\times}}
\def\wnto{{\widehat{\nto}}}

\def\nbo{{\zmin\otimes_\B \zmin}}
\def\nboplus{{\zmin^+\otimes_\B \zmin^+}}
\def\arith{{(\wnt,\zmax)}}
\def\wntb{{\widehat{\N^{\times 2}}}}
\def\arithb{{(\wntb,\nbo)}}
\def\arithc{{(\wntb,\conv(\Z\times \Z))}}
\def\beps{{\cR_\epsilon}}
\def\germ{{{\rm Germ}_{\epsilon=0}(\rmax)}}

\def\hatz{{\hat\Z^\times}}
\begin{document}

\title{Geometry of the  Arithmetic   Site}
\author{Alain Connes}
\email{alain@connes.org}
\address{Coll\`ege de France,
3 rue d'Ulm, Paris F-75005 France\newline
I.H.E.S. and Ohio State University.}
\author{Caterina Consani}
\email{kc@math.jhu.edu}
\address{Department of Mathematics, The Johns Hopkins
University\newline Baltimore, MD 21218 USA.}
%
%\dedication{}
\thanks{The second author would like to thank the Coll\`ege de France for hospitality and financial support.
}
%This is the abstract\newline (see \textbf{elsdoc.pdf} on the svn and/or one of the sample %files
%\code{cmguide1.pdf} or \code{cmguide1.ps} available at
%\url{http://www.compositio.nl/cmauthor.html}).

\begin{abstract}
We introduce  the {\em Arithmetic Site}: an algebraic geometric space deeply related to the non-commutative geometric approach to the Riemann Hypothesis.
We prove that the non-commutative space quotient of the ad\` ele class space of the field of  rational numbers by the maximal compact subgroup of the id\`ele class group,
which we had previously shown to yield the correct counting function to obtain the complete Riemann zeta function as Hasse-Weil zeta function, is  the set of geometric points of the arithmetic site over the  semifield of tropical real numbers. The action of the multiplicative group of positive real numbers on the ad\`ele class space corresponds to the action of the Frobenius automorphisms  on the above geometric points.  The underlying topological space of the arithmetic site is the  topos  of functors from the multiplicative semigroup  of non-zero natural numbers to the category  of sets. The  structure sheaf is made by semirings of characteristic one and is given globally  by the semifield of tropical integers.
 In spite of the countable  combinatorial nature of the arithmetic site, this space admits a one parameter semigroup  of Frobenius correspondences obtained  as sub-varieties of the square of the site.  This square   is  a semi-ringed topos whose structure sheaf involves Newton polygons. Finally, we show that the arithmetic site is intimately related to the structure of the (absolute) point in non-commutative geometry.\end{abstract}

\maketitle
\vspace*{1pt}
\tableofcontents  % for this guide only.
% A table of contents should normally not be included

\section{Introduction}

  It has long been known since \cite{Co-zeta} that the noncommutative space of ad\`ele classes of a global field provides a framework to interpret the explicit formulas of Riemann-Weil in number theory as a trace formula. In \cite{CC1}, we showed that if one divides the ad\`ele class space $\A_\Q/\Q^\times$ of the rational numbers by the maximal compact subgroup $\hatz$ of the id\`ele class group,  one obtains by considering the induced action of $\R_+^\times$, the counting distribution $N(u)$, $u\in [1,\infty)$, which determines, using the Hasse-Weil formula in the limit $q\to 1$, the complete Riemann zeta function. This analytic construction provides the starting point of the noncommutative attack to the Riemann Hypothesis. In order to adapt the  geometric proof of A. Weil, what was still missing until now was the definition of a geometric space of classical type whose points (defined over a ``field" replacing the algebraic closure of the  finite field $\F_q$ as $q\to 1$) would coincide with the afore mentioned quotient space. The expectation being that the action of suitably defined Frobenius automorphisms on these points would correspond to the above action of $\R_+^\times$.\newline
The primary intent of this paper is to provide a natural solution to this  search by introducing (\cf ~Definition \ref{site})  the {\em arithmetic site}  as an object  of algebraic geometry involving two elaborate mathematical concepts: the notion of topos and of (structures of) characteristic $1$ in algebra. 
 The topological space underlying the arithmetic site is the Grothendieck topos of sets with an action of  the multiplicative mono\"{\i}d $\nt$ of non-zero positive integers. The  structure sheaf of the arithmetic site is a 
 fundamental semiring of characteristic $1$, \ie $\Z_{\rm max}:= (\Z\cup\{-\infty\},\text{max},+)$ on which $\nt$ acts by  Frobenius endomorphisms. The role of the algebraic closure  of $\F_q$, in the limit $q\to 1$, is provided by the semifield $\rmax$ of tropical real numbers which  is endowed with a one parameter group of  Frobenius
automorphisms $\fr_\lambda$, $\lambda\in \R_+^\times$, given by $\fr_\lambda(x)=x^\lambda$ $\forall x\in \rmax$.\newline 
In this article we prove the following
\begin{thm} \label{structure3z} The set of points of the arithmetic site  $\arith$ over $\rmax$ coincides with the quotient of $\A_\Q/\Q^\times$ by the action of $\hatz$. The action of the Frobenius
automorphisms $\fr_\lambda$ of $\rmax$ on these points corresponds to the action of the id\`ele class group on $\hatz\backslash\A_\Q/\Q^\times$.
\end{thm}
The definition of the arithmetic site arises as a natural development of our recent work which underlined  the following facts 
 
$-$~The theory of topoi of Grothendieck provides the best geometric framework to understand cyclic (co)homology and the $\lambda$-operations using the (presheaf) topos associated to the cyclic category \cite{CoExt} and its epicyclic refinement (\cf \cite{topos1}).

$-$~Both the cyclic and the epicyclic categories, as well as the points of the associated topoi are best described from projective geometry over algebraic extensions of $\Z_{\rm max}$ (\cf\cite{topos}).

The arithmetic site acquires its algebraic structure from its structure sheaf.  In Section \ref{sectarithsite}, we describe the stalks of this sheaf at the points of the topos $\wnt$ which are in turn described in Section \ref{sectlb}, in terms of rank one ordered abelian groups. In Section \ref{secthw}, we combine the above theorem with our previous results as in \cite{CC1, CC2} to obtain  (\cf~Theorem \ref{mainthm}) a description of the complete Riemann zeta function $\zeta_\Q(s)=\pi^{-s/2}\Gamma(s/2)\zeta(s)$ as the Hasse-Weil zeta function of the arithmetic site. 
\begin{figure}
\begin{center}
\includegraphics[scale=.8]{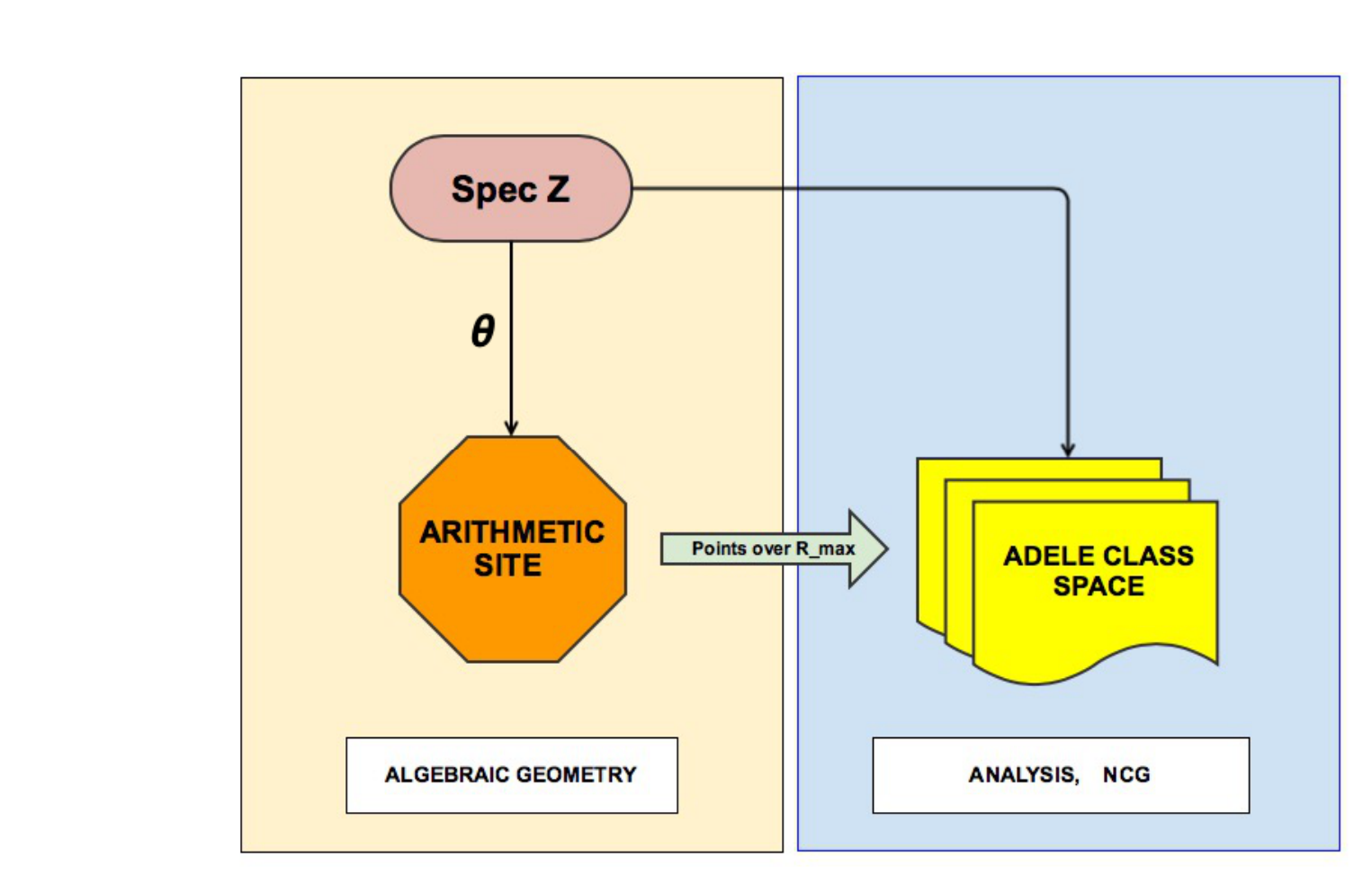}
\end{center}
\caption{The arithmetic site, $\Spec\Z$, and the adele class space.\label{clock1} }
\end{figure}
This construction makes heavy use of analysis and is naturally exploited on the right hand side of the correspondence between geometry and analysis described in Figure \ref{clock1}. The role of the prime numbers, as parameters for the periodic orbits of the Frobenius flow is refined in Section \ref{sectspz0} to construct (\cf~Theorem \ref{thmspz}) a geometric morphism of topoi between $\Spec\Z$ and $\wnto$. The slight difference between $\wnt$ and  the topos $\wnto$, dual to the semigroup $\nto$ of non-negative integers, is given by adjoining a base point to the points of $\wnt$: \cf~\S \ref{sectnto}. In Proposition \ref{sheavesonspz} we provide an interpretation of the pullback of the structure sheaf of the arithmetic site in terms of Cartier divisors on the spectrum of the ring of cyclotomic integers.

The general strategy adopted in this paper is to take full advantage of the knowledge achieved on both sides of the correspondence of Figure \ref{clock1}. The left hand side of that picture develops naturally into the investigation of the square of the arithmetic site and into the definition of the Frobenius correspondences. 
 In Section  \ref{sectsquare}, we describe the square of the arithmetic site as the topos $\wntb$ endowed with the structure sheaf defined globally by the tensor square $\nbo$ over the  smallest Boolean semifield of characteristic one. The semiring $\zmin:= (\Z\cup\{\infty\},\text{min},+)$ is isomorphic to $\zmax$ (by $n\mapsto -n$) but is more convenient for drawing figures. The reduced square is then obtained by reducing the involved semirings to become multiplicatively cancellative. In this way one achieves the important result of working with semirings whose elements are Newton polygons and whose operations are given by the convex hull of the union and the sum. In Proposition \ref{tensjustifsquare}, we prove that the points of the  square of the arithmetic site over $\rmax$ coincide with the product of the points of the  arithmetic site over $\rmax$. Then, we describe the Frobenius correspondences  as congruences on the square parametrized by positive real numbers $\lambda\in \R_+^\times$. The remarkable fact to notice at this point is that while the arithmetic site is constructed as a combinatorial object of countable nature it  possesses nonetheless a one parameter semigroup of ``correspondences" which can be viewed as congruences in the square of the site. In the context of semirings, the congruences \ie the equivalence relations compatible with addition and product, play the role of the ideals in ring theory. The Frobenius correspondences $\Psi(\lambda)$, for a rational value of $\lambda$, are deduced from the diagonal of the square, which is described by the product structure of the semiring, by  composition with the Frobenius endomorphisms. We interpret these correspondences geometrically, in terms of the congruence relation on Newton polygons corresponding to their belonging to the same half planes with rational slope $\lambda$. These congruences continue to make sense also for irrational values of $\lambda$ and are described using the best rational approximations of $\lambda$, while different values  of the parameter give rise to distinct congruences. In Section \ref{sectcomp}, we compute the composition law of the Frobenius correspondences  and we show the following
\begin{thm}\label{thmcomp0} Let $\lambda, \lambda' \in \R_+^\times$ be such that $\lambda\lambda'\notin \Q$. The composition of the Frobenius correspondences is then given by the rule
$$
\Psi(\lambda)\circ \Psi(\lambda')=\Psi(\lambda\lambda').
$$
The same equality still holds if $\lambda$ and $\lambda'$ are rational numbers. When  $\lambda, \lambda'$ are irrational and $\lambda\lambda'\in \Q$ one has
$$
\Psi(\lambda)\circ \Psi(\lambda')=\id_\epsilon\circ \Psi(\lambda\lambda')
$$
where $\id_\epsilon$ is the tangential deformation of the identity correspondence.
\end{thm}

Finally, in Section \ref{sectncg} we establish the link between the structure of the (absolute) point in noncommutative geometry and the topos $\wnt$. In particular, we recast the classification of matro\"ids obtained by J. Dixmier in \cite{Dixmier} in terms of the  non-commutative space  of points of the topos $\wnt$.
\subsection{Notations}

\subsubsection*{Topos}
The main reference for the theory of (Grothendieck) topoi is \cite{AGV}. We shall denote by $\pt$ the ``point" in topos theory \ie the topos of sets. Throughout the paper we also use extensively \cite{MM}. 

\subsubsection*{Characteristic one}
We denote by  $\rmax$ the semifield that plays a central role in idempotent analysis (\cf \cite{Maslov, Litvinov}) and tropical geometry (\cf \eg \cite{Gathmann, Mikhalkin}). It is the locally compact space $\R_+=[0,\infty)$ endowed with the ordinary product and the idempotent addition $x +' y=\max\{ x ,y \}$. This structure admits a one parameter group of automorphisms $\fr_\lambda\in \Aut(\rmax)$, $\fr_\lambda(x)=x^\lambda$ $\forall x\in\rmax$ which is the analogue of the arithmetic Frobenius in positive characteristic. The fixed points of the operator $\fr_\lambda$ form the Boolean semifield $\B=\{0,1\}\subset \rmax$: this is the only finite semifield which is not a field. One has 
${\rm Gal}_\B(\rmax)=\R^\times_+$.
In this paper  $\rmax$ denotes in fact the multiplicative version of the tropical semifield of real numbers $\R_{\rm max}:= (\R\cup\{-\infty\},\text{max},+)$.\newline
Let $\Z_{\rm max}:= (\Z\cup\{-\infty\},\text{max},+)$ be the  semifield of  tropical integers. Notice that $\Z_{\rm max}\simeq\zmin:= (\Z\cup\{\infty\},\text{min},+)$: the isomorphism mapping $\zmin\ni n\mapsto -n\in \Z_{\rm max}$.\newline
In this article, we shall use the multiplicative notation to refer to elements in $\zmin$, \ie we associate to $n\in\zmin$ the {\em symbol} $q^n$. In this way, the second operation of $\zmin$  becomes the ordinary product. If one represents $q$ as a positive, real number $0<q<1$, the first operation corresponds to the addition in $\rmax$: $x\vee y:=\max(x,y)$.\newline   
All these semirings are of characteristic $1$ \ie the multiplicative unit $1$ is idempotent for the addition and fulfills the equation $1\vee 1=1$. 
The operators $\ff_k\in \End(\zmin)$:  ~
%\begin{equation}\label{Fr}
$\N^\times \to \End(\zmin),~k\mapsto\ff_k(n) := kn$ 
%\end{equation}
are the analogues, in characteristic $1$, of the Frobenius endomorphism in characteristic $p>1$ and, in the multiplicative notation, they are defined by the rule $\ff_k(x)=x^k$.

\section{The points of the topos $\wnt$}\label{sectlb}

In this section we show that the category of  points of the topos $\wnt$ is equivalent to the category of rank one ordered groups, \ie of totally ordered groups isomorphic to non-trivial subgroups of   $(\Q,\Q_+)$, and injective morphisms of ordered groups. We recall that for a topos of presheaf type $\hat\cC$, \ie the topos of contravariant functors $G:\cC\longrightarrow\Se$ from a small category $\cC$ to the category $\Se$ of sets, any object $C$ of the category $\cC$ defines a point $p$ of $\hat\cC$ whose pullback  $p^*$ is given  by the evaluation  $G\mapsto p^*(G):= G(C)$. Moreover, it is well known that every point of $\hat\cC$ is obtained as a filtering colimit of points of the above form. In particular, for the topos $\wnt$ this implies that every point can be obtained from  a (filtering) sequence $(n_i)$, with $n_i\in\nt$ and $n_i\vert n_{i+1}$ $\forall i$. Two cofinal sequences label the same point. The equivalence   relation $(n_i)\sim (m_i)$ on  sequences  to define isomorphic points states that the classes of the limits 
$n=\lim n_i, ~m=\lim m_i \in \hat\Z/\hatz$ are the same in the double quotient $\Q_+^\times\backslash\A^f/\hat\Z^\times
$, where $\A^f= \hat\Z\otimes_\Z\Q$ denotes the finite ad\`eles of $\Q$ and where $\Q_+^\times$ acts by multiplication on $\A^f$. In this section we shall give  a detailed account of this construction of the points of the topos $\wnt$ and discuss, also, the case of the topos $\wnto$.

  \subsection{The points of a presheaf topos} \label{preshtop}
 It is a standard fact in topos theory that the category of points a topos  $\hat \cC$, where $\cC$ is a small category,  is canonically equivalent to the category of {\em flat} functors $\cC\longrightarrow \Se$ and natural equivalences between them. 
We recall that a covariant functor $F:\cC\longrightarrow \Se$ is flat if and only if it is filtering \ie $F$ satisfies the following conditions
\begin{enumerate}
\item $F(C)\neq \emptyset$, for at least one object $C$ of $\cC$.
\item Given two objects $A,B$ of $\cC$ and elements $a\in F(A)$, $b\in F(B)$, there exists
an object $Z$ of $\cC$, morphisms $u:Z\to A$, $v:Z\to B$ and an element $z\in F(Z)$, such that 
$
F(u)z=a, \   F(v)z=b.
$
\item Given two objects $A,B$ of $\cC$ and arrows $u,v:A\to B$ and $a\in F(A)$ with $F(u)a=F(v)a$, there
exists an object $Z$ of $\cC$, a morphism $w:Z\to A$,  and an element $z\in F(Z)$, such that 
$$
F(w)z=a, \  \  u\circ w=v\circ w\in \Hom_\cC(Z,B).
$$
\end{enumerate}

\subsection{Points of $\wnt$ and ordered groups}\label{ptstogr}
We shall investigate the case of the small category $\cC = \N^\times$ with a single object $\ast$. A point of $\wnt$ is described by a covariant functor $F: \N^\times\longrightarrow \Se$ which is filtering, \ie such that the category $\int_{\N^\times}\, F$ is filtering (\cf \cite{MM} Chapter VII \S 6 for the notation). The functor $F$ is determined by the set $X=F(*)$ endowed with an action of the semi-group $\N^\times$: $F(k):X\to X$, $\forall k>0$. The category $\int_{\N^\times}\, F$ has objects given by elements $x\in X$ and the morphisms between any pair of objects $x,y\in X$ are provided by the elements of $\nt$ such that $F(k)x=y$. The filtering condition for the category $\int_{\N^\times}\, F$  means that
\begin{enumerate}
\item $X\neq \emptyset$.
\item For any $x,x'\in X$ there exist $z\in X$ and $k,k'\in \N^\times$ such that $F(k)z=x$, $F(k')z=x'$.
\item For $x\in X$ and $k,k'\in \N^\times$, the equality $F(k)x=F(k')x$ implies $k=k'$.
\end{enumerate}
The last condition here above follows from the third condition of \S \ref{preshtop} and the fact that the semigroup $\nt$ is simplifiable.

Next theorem provides an algebraic description of the points of the topos $\widehat{ \N^\times} $.

\begin{thm}\label{thmstart}
The category of points of the topos $\widehat{ \N^\times} $ is canonically equivalent to the category of totally ordered groups isomorphic to non-trivial subgroups of   $(\Q,\Q_+)$ and injective morphisms of ordered groups.
\end{thm}
\proof The proof of this theorem follows from the next two lemmas and some final, easy considerations.

\begin{lem}\label{flat1} Let $F: \N^\times\longrightarrow \Se$ be a flat functor. The following equality defines a commutative and associative addition  on $X=F(\ast)$ 
\begin{equation}\label{defnsum}
x+x':=F(k+k')z, \ \forall z\in X, \  F(k)z=x, \  F(k')z=x'.
\end{equation}
Endowed with this operation $X$ coincides with the strictly positive part of an abelian, totally ordered group  $(H,H_+)$ which is an increasing union of subgroups $(\Z,\Z_+)$.
\end{lem}
\proof We show that the operation $+$ on $X$ is well defined, that means independent of the choices. 
If for some $y\in X$ and $\ell,\ell'\in \nt$ one has  $F(\ell)y=x$, $F(\ell')y=x'$, one uses $(ii)$ of the filtering condition to find $u\in X$ and $a,b\in \N^\times$ such that 
$
F(a)u=y, \, F(b)u=z
$. 
One then has 
$
F(\ell a)u=x=F(kb)u, ~ F(\ell' a)u=x'=F(k'b)u
$
and then one uses $(iii)$ of the filtering condition to obtain 
$
a\ell=bk,  ~a\ell'=bk'
$.
Thus one gets
$$F(k+k')z=F((k+k')b)u= F(kb+k'b)u=F(a\ell+a \ell')u=F(\ell+\ell')F(a)u=F(\ell+\ell')y.$$
This shows that the addition is well defined. This operation is commutative by construction. Notice that for any finite subset $Z\subset X$ one can find $z(Z)\in X$ such that $Z\subset F(\nt) z(Z)$. The associativity of the addition then follows from the associativity of the addition of integers. The obtained additive semigroup $(X,+)$ is therefore an increasing union
$$
(X,+)=\displaystyle{\bigcup_{Z\subset X\atop Z \rm{finite}}} z(Z)[1,\infty)
$$ 
of the additive semigroups $z(Z)[1,\infty)$, each isomorphic to the  additive semigroup of integers $n\geq 1$.
From this fact one derives that the obtained semigroup is simplifiable and that for any  pair  $a,b\in X$, $a\neq b$, there  exists $c\in X$ such that $a+c=b$ or $b+c=a$. By symmetrization of  $(X,+)$ one obtains  an abelian, totally ordered group  $(H,H_+)$ which is an increasing union of subgroups isomorphic to $(\Z,\Z_+)$. One has $X=\{h\in H\mid h>0\}$. \endproof

Notice that the additive structure $(X,+)$ determines the action of  $\N^\times$ on $X=F(\ast)$ by the equation
$$
kx= x+\cdots +x \  \  (k \ \text{terms}).
$$

\begin{lem}\label{flat2} The ordered abelian group $(H,H_+)$ obtained as the symmetrization of $(X,+)$ is isomorphic to a subgroup of $(\Q,\Q_+)$. For any two injective morphisms  $j, j': H\to  \Q$, there exists $r\in \Q^\times_+$ such that $j'=rj$.
\end{lem}
\proof
Let $x\in X$: we show that there exists a unique injective morphism $j_x: H\to \Q$  such that $j_x(x)=1$. For any $x'\in X$ and with the notations of \eqref{defnsum} one sets 
$$
j_x(x'):=k'/k, \ \forall z\in X, \  F(k)z=x, \  F(k')z=x'.
$$
This is well defined since, with the notations of the proof of Lemma \ref{flat1}, another choice $y,\ell,\ell'$ of $z,k,k'$ 
gives the equalities $a\ell=bk$, $a\ell'=bk'$ showing that $k'/k=\ell'/\ell$ is independent of any choice. One easily checks that $j_x$ is also additive and injective. \newline 
We claim that given two subgroups $H, H'\subset \Q$, any non-trivial ordered group morphism $\phi:H\to H'$ is of the form $\phi(x)=rx$, $\forall x\in H$ for some $r\in \Q^\times_+$. Indeed, one takes $x_0\in H$, $x_0>0$, and lets $r\in \Q^\times_+$ such that $\phi(x_0)=rx_0$. Then, for any $x\in H_+$ there exist integers $n>0$ and $m\geq 0$ such that $nx=mx_0$. Then one derives $n\phi(x)=\phi(nx)=\phi(mx_0)=m\phi(x_0)=mrx_0=nrx$ and thus $\phi(x)=rx$. \endproof

To finish the proof of the theorem, one notices that a natural transformation of flat functors is by definition an $\N^\times$-equivariant map $f:X\to X'$. Then, the properties $(ii)$ and $(iii)$ of flat functors show that $f$ is necessarily injective. Moreover,  for $x,y\in X$, $x=F(n)z$, $y=F(m)z$, one has  
$$
f(x+y)=f(F(n+m)z)=F'(n+m)f(z)=f(x)+f(y)
$$
using $F'(n)f(z)=f(x)$ and $F'(m)f(z)=f(y)$.
Thus, $f$ is an injective morphism of ordered groups. Conversely, an injective morphism of ordered groups gives a natural transformation of the associated flat functors. Thus one finally concludes that the category of points of the topos $\widehat{ \N^\times} $ is the category of totally ordered abelian groups of rank one (isomorphic to subgroups of $\Q$) and injective morphisms.
\endproof

\subsection{The case of $\nto$}\label{sectnto}

We investigate a variant of  Theorem \ref{thmstart} obtained by replacing the multiplicative mono\"{\i}d $\nt$ with the (pointed) mono\"{\i}d $\nto$ defined by adjoining a $0$-element: $\nto$ is the multiplicative mono\"{\i}d of non-negative integers. Let $\cC'$ be the  small category  with one object $\ast$ and endomorphisms given by $\nto$.  As before, a flat functor $F:\cC'\longrightarrow \Se$ is described by assigning a set $X=F(\ast)$ and an action $F$ of the mono\"{\i}d $\nto$ on $X$ fulfilling the filtering conditions. Let $P=F(0)$, this determines a map $P:X\to X$ which fulfills
$$
P\circ P=P, \   \  P\circ F(n)=F(n)\circ P=P\qqq n\in \nt.
$$
Thus, the image of $P$ is described by the subset $X_0=\{x\in X\mid Px=x\}$. We use the filtering conditions to investigate $X_0$. The second  condition states that
for any pair $x,x'\in X$ there exist $z\in X$ and $k,k'\in \nto$ such that $F(k)z=x$, $F(k')z=x'$. If   both $x,x'\in X_0$, then one gets
$$
x=Px=PF(k)z=Pz, \  \ x'=Px'=PF(k')z=Pz
$$
and  hence $x=x'$. Moreover the first filtering condition shows that the set $X$ is non-empty and  so is $X_0$,  thus it contains exactly one element that we denote $0_X$. Next, let $ x\in X$, $n\in \nt$ and assume that for some $n\in \nt$ one has $F(n)x=0_X$. Then one has $F(n)x=F(0)x$ and this time the third filtering condition states the existence of a morphism $w:*\to *$,  and an element $z\in X$, such that 
$
F(w)z=x, ~n\circ w=0\circ w\in \Hom_\cC(*,*)
$. Since $n\in \nt$ one has 
 $w=0$, and the equation $F(w)z=x$ implies that $x=0_X$.  It follows that the set  $X^*:=X\setminus \{0_X\}$ is stable under the action of $\nt$. Let us assume that $X^*\neq \emptyset$. We show that  the restriction of the action of $\nt$ on $X^*$ fulfills   the filtering conditions. The second filtering condition is fulfilled since given two  elements $a,b\in X^*$, there exists
 morphisms $u, v\in \nto$ and an element $z\in X$, such that 
$
F(u)z=a, ~  F(v)z=b.
$
The element $z\in X$ cannot be $0_X$ since $F(n)0_X=0_X$ for all $n\in \nto$. Thus $z\in X^*$. Moreover one cannot have $u=0$ or $v=0$ since we assume $a,b\in X^*$. The third filtering condition for the flat functor $F:\cC'\longrightarrow \Se$ states that given   $u,v\in \nto$ and $a\in X$ with $F(u)a=F(v)a$ there
exists  $w\in \nto$,  and an element $z\in X$, such that 
$
F(w)z=a, ~  u\circ w=v\circ w\in \nto
$.
Let us assume that $a\in X^*$ and that $u,v\in \nt$. Then the condition $F(w)z=a$ shows that one cannot have $w=0$ since this would imply that $a=0_X$. Thus one has $w\in \nt$ and $u=v$. This shows that the restriction of the action of $\nt$ to $X^*$ fulfills the three filtering conditions, provided one assumes  $X^*\neq \emptyset$. 

With this verification, the replacement for Theorem \ref{thmstart} is given by the following
\begin{thm}\label{thmstartbis}
The category of points of the topos $\widehat{ \nto} $ is canonically equivalent to the category of totally ordered groups isomorphic to  subgroups of   $(\Q,\Q_+)$, and  morphisms of ordered groups.
\end{thm}
\proof Given a point of $\wnt$ \ie a flat functor $F:\nt\longrightarrow \Se$, we consider the set $\tilde X:=X\cup \{0_X\}$. We extend the action of $\nt$ on $X$ to an action of $\nto$ on $\tilde X$ such that $0_X$ is fixed by all $F(n)$ and setting $F(0)y=0_X$, $\forall y\in X$. This extended functor $\tilde F$ fulfills the first two filtering conditions. We check that the third also holds. Let $u,v\in \nto$ and $a\in \tilde X$ with $F(u)a=F(v)a$. Then, if $a=0_X$ we take the morphism $w=0\in \nto$ and $z=a$. One has $F(w)z=a$ and $u\circ w=v\circ w\in\nto$. If $a\neq 0_X$ and $F(u)a=F(v)a=0_X$, it follows that necessarily one has $u=v=0\in \nto$. Thus, one can take $w=1$ and $z=a$. Finally, if $a\neq 0_X$ and $F(u)a=F(v)a\neq 0_X$, then one has $u,v\in \nt$ and one uses the flatness for the action of $\nt$ on $X$. This proves that $\tilde F$ is flat and hence defines a point of $\wnto$. Moreover the above discussion shows that any point of $\wnto$ different from the one element set $X=\{0_X\}$ is obtained in this manner. We claim that under an equivariant map $f:X\to X'$ of $\nto$-sets one has $f(0_X)=0_{X'}$ and if $f(x)=0_{X'}$ for some $x\in X$, $x\neq 0_X$ then one has $f(y)=0_{X'}$ for all $y\in X$. Indeed, there exists $z\in X$ and $n,m$ such that $F(n)z=x$, $F(m)z=y$. One has $n\in \nt$ since $x\neq 0_X$ and thus one gets $f(z)=0_{X'}$ since otherwise $f(x)\neq 0_{X'}$ as the complement of $0_{X'}$ is stable under $F'(n)$ for $n\in \nt$. This shows that the only new morphism is the morphism $0$ and the only new point is the point $0$. \qed\vspace{.1in}

Notice that the category of points of the topos $\widehat{ \nto} $ is {\em pointed}, \ie it admits a (unique) initial and final object $0$. The above proof shows  that the category of points of the topos $\widehat{ \nto} $ is  obtained canonically by adjoining to the category of points of $\wnt$ an object which is both initial and final.
\subsection{Adelic interpretation}\label{sectadel}

The ring  $\A^f$ of finite ad\`eles of $\Q$ is defined as the restricted product $\A^f=\prod' \Q_p$: its maximal compact subring is the profinite completion $\hat\Z=\prod \Z_p$ of $\Z$.  One has $\A^f=\hat\Z\otimes_\Z\Q$ and there is a canonical embedding $\Q\subset \A^f$ such that $\Q\cap \hat \Z=\Z$. Let $\Q_+^\times$ act by multiplication on $\A^f$.

The following result will play an important role in the paper.

\begin{prop}\label{propsubgroup} $(i)$~Any non-trivial subgroup of $\Q$ is uniquely of the form
\begin{equation}\label{param}
H_a:=\{q\in \Q\mid aq\in \hat\Z\}, \  \ a\in \A^f/\hatz
\end{equation} 
where $\hatz$ denotes the multiplicative group of invertible elements in the ring $\hat \Z$ acting by multiplication on $\A^f$.

$(ii)$~The map $a\mapsto H_a$ of \eqref{param} induces a canonical bijection of  the quotient space
\begin{equation}\label{bc}
\Q_+^\times\backslash\A^f/\hatz
\end{equation}
  with the set of isomorphism classes of points of the topos $\wnt$.
\end{prop}
\proof $(i)$~We give a conceptual proof using Pontrjagin duality for abelian groups. \newline
 Let $\chi:\Q/\Z\to U(1)$, $\chi(\alpha)=e^{2\pi i \alpha}$. It  gives a character, still noted $\chi$ of $\A^f$, whose kernel is the additive subgroup $\hat\Z\subset \A^f$ (use the canonical isomorphism $\Q/\Z=\A_f/\hat \Z$). Then the  pairing 
$$
<q,a>=\chi(qa)\qqq q\in \Q/\Z, \  a\in \hat\Z
$$
 identifies $\hat\Z$ with the  Pontrjagin dual of the discrete abelian group $ \Q/\Z$. To prove $(i)$ we can assume that the non-trivial subgroup of $\Q$ contains $\Z$. The subgroups of $
\Q$ which contain $\Z$ are in fact determined by subgroups $H\subset \Q/\Z$. Given a subgroup $H\subset \Q/\Z$ one has
$$
H=(H^\perp)^\perp=\{q\in \Q/\Z\mid  <q,x>=1,\  \forall x\in H^\perp\}, \ \ H^\perp=\{x\in \hat\Z\mid <q,x>=1,\  \forall q\in H\}.
$$
Let $J=H^\perp$. It is a closed subgroup $J\subset \hat\Z$. Then, the equality $H=H_a$ for some $a\in \hat \Z$ derives from the following observations on closed subgroups $J\subset \hat\Z$, by showing that they are all of the form $J=a\hat\Z$, for some $a\in \hat \Z$ unique up to the action of $\hatz$ by multiplication.
\begin{itemize}
\item A closed subgroup of $\hat\Z$ is an ideal in the ring $\hat\Z$.
\item For each prime $p$, the projection $\pi_p(J)\subset \Z_p$ coincides with the intersection $(\{0\}\times \Z_p)\cap J$
and is a closed ideal $J_p\subset\Z_p$, moreover one  has: $x\in J\iff \pi_p(x)\in J_p$ $\forall p$.
\item Any ideal of $\Z_p$ is principal.
\end{itemize}
Thus one writes $J_p=a_p\Z_p$ for all $p$ and  it follows that $J=a\hat\Z$ where $a=(a_p)\in \prod\Z_p$. Moreover, for each $p$ the element $a_p\in\Z_p$ is unique up to the action of $\Z_p^\times$ and in this way one also proves the statement on the uniqueness.\newline
$(ii)$~By Theorem \ref{thmstart} any point of the topos $\wnt$ is obtained from a subgroup $H_a\subset \Q$. Elements $a,b\in \A^f/\hatz$ determine isomorphic points of the topos $\widehat{ \N^\times} $
if and only if the ordered groups $H_a$ and $H_b$ are abstractly isomorphic. An abstract isomorphism is given by the multiplication by  $q\in \Q_+^\times$ so that $H_b=qH_a$. It follows from $(i)$ that $a=qb$ in $\A^f/\hatz$ and one derives the result. \endproof

\begin{rem}\label{natnum}{\rm 

The quotient of the additive group $\A^f$ by the open compact subgroup $\hat\Z$ is the direct sum
$$
\Q/\Z=\A^f/\hat \Z= \bigoplus \Q_p/\Z_p.
$$
When this fact is applied to a rational number $q\in \Q$, this yields the decomposition in simple elements
\begin{equation}\label{decsimple}
q=n+
\sum  n_{p_j}/p_j^{\alpha_j}
\end{equation} 
where $n\in \Z$, the $p_j$'s are prime numbers, $0<n_{p_j}<p_j^{\alpha_j}$, and $(n_{p_j},p_j)=1$. One can use this decomposition to give a direct proof of Proposition \ref{propsubgroup}.
}\end{rem}

\begin{rem}\label{natnum1}{\rm 

One can describe the elements of $\hat \Z/\hat \Z^\times$ (which parametrize the subgroups $\Z\subset H\subset \Q$) by means of ``supernatural numbers" in the sense of Steinitz. This is the point of view adopted in \cite{Dixmier, LB}. By definition, a supernatural number is a formal product $\prod p^{n_p}$ over all primes, where the exponents $n_p$ belong to $\N\cup \{\infty\}$. One can check easily that this product is convergent in $\hat \Z/\hat \Z^\times$ and that one obtains in this way a canonical bijection of the supernatural numbers with $\hat \Z/\hat \Z^\times$ compatible with the labeling of subgroups $\Z\subset H\subset \Q$.
}\end{rem}

\begin{rem}\label{permsect}{\rm 

In  \S \ref{sectarithsite}, we shall enrich the topos  $\widehat{ \N^\times} $ with the structure of a curve over the Boolean semifield $\B$, by  endowing the topos with its structure sheaf.  The presence of a large group of automorphisms on the topos $\wnt$, arising from the automorphisms of the mono\"{i}d $\nt$, shows that one does not have, at this state of the construction, a sufficient geometric structure on the topos. It is however useful to investigate how the endomorphisms $\rho$ of the mono\" id $ \N^\times$ act on the  points of  $\widehat{ \N^\times}$ in terms of the above adelic description.
Any such  $\rho$ extends by continuity to  an endomorphism $\bar \rho$ of the multiplicative mono\" id $\hat \Z/\hatz$ and this extension determines the action of  $\rho$ on the  points of  $\wnt$.
The action of $\rho$ on the point associated  to the ordered group  $H=H_a$,  $a\in \hat \Z/\hatz$ is given by the point associated  to the ordered group  $H=H_b$, with  $b=\bar \rho(a)\in \hat \Z/\hatz$.
}\end{rem}

\begin{rem}\label{topspec}{\rm In \S \ref{sectspz0}, we shall construct a geometric morphism of topoi $\Theta:\Spec\Z\to\wnto$. When considered as a map of points, it associates to the point of $\Spec\Z$ corresponding to 
 a prime $p$ the point of $\wnto$ associated to the class of  the finite adele $\alpha_p=p^\infty\in\A^f/\hatz$  whose components are all equal to $1$ except at $p$ where the component vanishes. One can show that the topology induced on these points by the quotient topology of $\Q_+^\times\backslash\A^f/\hatz$ coincides with the Zariski topology of $\Spec\Z$, \ie the open sets are the complements of the finite subsets. The geometric morphism $\Theta$ associates to the generic point of $\Spec\Z$ the base  point of $\wnto$.
}\end{rem}

\section{The Arithmetic Site $\arith$}\label{sectarithsite}
In this section we introduce  the structure  sheaf $\zmax$ on the topos $\wnt$. This additional structure on  $\wnt$, without which the group of automorphisms (\cf~Remark \ref{permsect}) would contain arbitrary permutations of the primes, turns the topos $\wnt$ into the arithmetic site $\arith$.

\subsection{The structure sheaf $\zmax$ and its stalks}
We start by introducing the definition of the arithmetic site.

\begin{defn}\label{site} The arithmetic site $\arith$ is the topos $\wnt$
endowed with the  {\em structure sheaf} $\cO:=\zmax$ viewed as a semiring in the topos using the action of $\nt$ by the Frobenius endomorphisms.
\end{defn}
The action of $\nt$ on $\zmax$ is by the Frobenius endomorphisms 
$\N^\times \to \End(\zmax),~k\mapsto\ff_k(n) := kn$.  
Likewise an algebraic scheme is a ringed space of a certain kind, the arithmetic site can be thought of as a ``semi-ringed topos".
The semiring structure of $\zmax$ is compatible with the action of $\nt$ and automatically endows the stalks   
of the sheaf $\zmax$ on the topos $\wnt$ with a particular structure of semiring which we now determine. 
\begin{thm} \label{structure2} The stalk of the structure sheaf $\cO$ at the point of the topos $\wnt$ associated to the ordered group $H\subset \Q$  is canonically isomorphic to the semiring 
$
H_{\text{max}}:=(H\cup\{-\infty\},\text{max},+)
$.
\end{thm}
\proof To the point of the topos $\wnt$ associated to the rank one ordered group $H$ corresponds the flat functor $F:\nt\longrightarrow\Se$ which associates to the single object $\ast$ of the small category $\nt$ the set $F(\ast)=H_{>0}$ and to the endomorphism $k$ the multiplication $F(k)$ by $k$ in $H_{>0}$. The inverse image functor connected with this point is the functor which associates to  any $\nt$-space $X$ its geometric realization  $|X|_F$ of the form
\begin{equation}\label{geometricreal}
   |X|_F=\left(X\times_{\nt} H_{>0} \right)/\sim.
\end{equation}
The relation $\sim$ states the equivalence of the pairs $(x,F(k)y)\sim (k.x,y)$, where $k.x$ denotes the action of $k$ on $x\in X$. More precisely, by using the flatness of the functor $F$ the equivalence relation is expressed as follows
$$
(x,y)\sim (x',y')\iff \exists z\in Y, k,k'\in \nt, \ F(k)z=y, \ F(k')z=y', \ k.x=k'.x'.
$$
 To obtain a description of the stalk of the structure sheaf $\cO$ at $H$, one applies \eqref{geometricreal} to the $\nt$-space $X=\zmax$, with the action of $\nt$ defined by $k.x:=kx$ and $k.(-\infty):=-\infty$. We define  
\begin{equation}\label{mapbeta}
\beta:|\zmax|_F\to H_{\rm max}, \  \ \beta(x,y)=xy\in H_{\rm max}\qqq x\in \zmax, y\in H_{>0}.
\end{equation}
The map $\beta$ is compatible, by construction, with the equivalence relation $\sim$. One has 
$$
\beta(1,y)=y, \ \beta(0,y)=0, \ \beta(-1,y)=-y, \ \beta(-\infty,y)=-\infty\qqq y\in H_{>0}
$$
which shows that $\beta$ is surjective. We show that $\beta$ is also injective. Let $u=\beta(x,y)=\beta(x',y')$, then if $u>0$ one has $x\in \nt$ and $(x,y)\sim (1,xy)=(1,x'y')\sim (x',y')$. For $u<0$, one has 
$-x\in \nt$ and $(x,y)\sim (-1,-xy)=(-1,-x'y')\sim (x',y')$. For $u=0$ one has $x=x'=0$ and one uses the flatness of $F$ to find $z\in H_{>0}$ and $k,k'\in \nt$ with $kz=y$, $k'z=y'$. Then one derives
$(0,y)= (0,kz)\sim (0,z)\sim (0,k'z)=(0,y')$. The same proof works also for $u=-\infty$.

The two operations $\vee$ and $+$ on $\zmax$  determine maps of $\nt$-spaces $\zmax\times \zmax\to\zmax$. Since the geometric realization functor commutes with finite limits one 
gets $|\zmax\times \zmax|_F=|\zmax|_F\times |\zmax|_F$. One  obtains corresponding maps 
$|\zmax|_F\times |\zmax|_F\to |\zmax|_F$ associated respectively to the two operations. The identification $|\zmax\times \zmax|_F=|\zmax|_F\times |\zmax|_F$ is described by the map
$$
(x,x',y)\mapsto (x,y) \times (x',y)\qqq x,x'\in\zmax, \  y\in H_{>0}.
$$
 One has $y(x\vee x')=yx\vee yx'$ for any $y\in H_{>0}$ and $x,x'\in\zmax$ and this suffices to show that the map $|\zmax|_F\times |\zmax|_F\to |\zmax|_F$ associated to the operation $\vee$ on $\zmax$ is the operation $\vee$ on $H_{\rm max}$. Similarly one has $y(x + x')=yx+ yx'$ for any $y\in H_{>0}$ and $x,x'\in\zmax$ and again this suffices to show that the map $|\zmax|_F\times |\zmax|_F\to |\zmax|_F$ associated to the operation $+$ on $\zmax$ is the operation $+$ on $H_{\rm max}$.
\endproof
\begin{rem}\label{caseno}{\rm The statement of Theorem \ref{structure2} continues to hold if one replaces the pair $\arith$ by $(\wnto,\zmax)$, where the action of $0\in \nto$ on $\zmax$ is the endomorphism $x\mapsto 0x=0$ $\forall x\in\Z\subset \zmax$ and $0(-\infty):=-\infty$. The map $\beta$ of \eqref{mapbeta} extends to the case $y\in H_{\geq 0}$ by setting $0(-\infty):=-\infty$. In the proof of injectivity of $\beta$, in the case $u=0$, one uses that $x$ or $y$ vanishes to get $(x,y)\sim (0,0)\sim (x',y')$. The proof in the other cases is unchanged. Moreover at the only new point of $\wnto$ \ie the point associated to the trivial group $\{0\}$, the stalk is the semifield $\B$.
}\end{rem}
\begin{rem}\label{casenote}{\rm  Theorem \ref{structure2} continues to hold if one replaces the pair $\arith$ by $(\wnt,\zminp)$, where $\zminp$ is the sub-semiring $\zminp:=\{n\in\zmin, n\geq 0\}$. The stalk at the point of the topos $\wnt$ associated to the ordered group $H\subset \Q$  is canonically isomorphic to the semiring 
$
R=(H_{\geq 0}\cup\{\infty\},\text{min},+)
$. This is the notation adopted in Theorem 2.4 of \cite{CCas}, where $\zminp$ was denoted $\overline\N$. One switches back and forth from the semiring $\zmax\sim \zmin$ to the semiring $\zminp\sim \overline\N$, by applying the following two functors. The first one, noted $\cO$, associates to a semiring $R$ of characteristic one the (sub) semiring $\cO(R):=\{x\in R\mid x\oplus 1=1\}$ where $\oplus$ denotes the addition in $R$ and $1$ is the multiplicative unit. The second functor, noted ${\rm Frac}$, associates to a semiring $R$ without zero divisors its semiring of fractions.
}\end{rem}

One should not confuse the semiring $\zmax$ with the semiring of  global sections of the structure sheaf of the arithmetic site. The latter selects the elements of $\zmax$ which are invariant under the action of $\nt$ as explained by the following
\begin{prop}\label{propglobsect} The global sections $\Gamma(\wnt,\zmax)$ of the structure sheaf  are given by the sub-semifield $\B\subset \zmax$.
\end{prop}
\proof For a Grothendieck topos $\cT$ the global section functor $\Gamma:\cT\to \Se$,  is given by $\Gamma(E):=\Hom_\cT(1,E)$, where $E$ is an object of $\cT$ and $1$ is the final object in the topos.
In the particular case of a topos of the form $\hat \cC$, where $\cC$ is a small category, a global section of a contravariant functor $P: \cC\longrightarrow\Se$ is a function which assigns to each object $C$ of $\cC$ an element $\gamma_C\in P(C)$ in such a way that for any morphism $f:D\to C$ one has $P(f)\gamma_C=\gamma_D$
(see \cite{MM} Chapter I \S 6 (9)). We apply this definition to $\cC=\nt$ and derive that the global sections coincide with the elements of $\zmax$ which are invariant under the action of $\nt$, \ie $\B\subset \zmax$.\endproof

Proposition \ref{propglobsect} indicates that, as a generalized scheme over $\B$, the arithmetic site is, inasmuch as it can be considered as the analogue of a curve, both {\em complete and irreducible}.

The above Theorem \ref{structure2} shows that one does not require additional structure in order to obtain the semiring structure on the stalks of the structure sheaf at the points of $\wnt$. In the next sections, we shall see that this structure has the right properties to obtain the Riemann zeta function and the counting function of the mysterious ``curve" underlying the geometry of the prime numbers.

\subsection{The points of the arithmetic site over $\rmax$}\label{sectrmax}

By definition (\!\!\cite{AGV} Expos\'e 4, Definition 13.1), a morphism of  ringed topos $(X,\cO_X)\to (Y,\cO_Y)$ is a pair $(f,\theta)$, where $f:X\to Y$ is a geometric morphism   and $\theta: f^*(\cO_Y)\to \cO_X$ is a morphism of rings on the topos $X$. Using the adjunction formula between the pullback $f^*$ and the push-forward $f_*$, one can equivalently describe $\theta$ as the morphism $\theta: \cO_Y\to f_*(\cO_X)$. This definition extends immediately to semi-ringed topoi, in this paper we apply it to the semi-ringed topos $(\pt,\rmax)$, where $\pt$ is the single point, \ie the topos of sets, endowed with the semifield $\rmax$.

\begin{defn} \label{defnpt} Let $(Y,\cO_Y)$ be a semi-ringed topos. A point of $(Y,\cO_Y)$ defined over $\rmax$ is a pair given by a point $p$ of $Y$ and a  morphism $f^\#_p: \cO_{Y,p}\to \rmax$ of semirings from the stalk of $\cO_Y$ at $p$ to $\rmax$.  
\end{defn}

Likewise a point of a topos $\cT$  is an equivalence class of isomorphic objects in the category of geometric morphisms $\pt\to \cT$, we say that two  points  of $(Y,\cO_Y)$ defined over $\rmax$  are equivalent $(p,f^\#_p)\sim (q,g^\#_q)$ when there exists an isomorphism of points $\phi:p\stackrel{\sim}{\to} q$ compatible with the morphisms of semirings \ie such that $g^\#_q\circ \phi^\#=f^\#_p$.

Before stating the main theorem of this section we introduce the following lemma which provides a canonical bijection between the quotient $\Q_+^\times\backslash((\A^f/\hatz)\times \R_+^\times)$ and the set of non-trivial rank one subgroups of $\R$.

\begin{lem}\label{subbis} Let $\Phi$ be the map from $(\A^f/\hatz)\times \R_+^\times$ to subgroups of $\R$ defined by 
$$
\Phi(a,\lambda):=\lambda H_a\qqq a\in \A^f/\hatz, \, \lambda \in  \R_+^\times.
$$
Then, $\Phi$ is a bijection between the quotient of $( \A^f/\hatz)\times \R_+^\times$ by the diagonal action of $\Q_+^\times$ and the set of non zero subgroups of $\R$ whose elements are pairwise commensurable.
\end{lem}
\proof The map $\Phi$ is invariant under the diagonal action of $\Q_+^\times$ since by replacing $(a,\lambda)$ with $ (qa,q\lambda)$ for $q\in \Q_+^\times$ one replaces $H_a$ by $H_{qa}=q^{-1}H_a$ and $\lambda H_a$ by 
$q\lambda q^{-1}H_a=\lambda H_a$. \newline
We show that $\Phi$ is surjective. 
A subgroup $H\subset \R$ whose elements are pairwise commensurable is, as an ordered group, isomorphic to a subgroup of $\Q$. Given $x\in H$, $x>0$  one can find  $\lambda\in  \R_+^\times$ such that $\lambda^{-1} x\in \Q$ and it follows that $\lambda^{-1} H\subset \Q$. This determines an ad\`ele  $a\in  \A^f/\hatz$ such that 
$ H=\lambda H_a$ and thus $H=\Phi(a,\lambda)$. \newline
Next, we prove that $\Phi$ is injective. Assume that $\lambda H_a=\lambda' H_b$, for $\lambda, \lambda' \in  \R_+^\times$. Then since both $H_a$ and $H_b$ are non-zero subgroups of $\Q$, one has $\lambda'=q\lambda$ for some $q\in \Q_+^\times$. One then gets $\lambda H_a=\lambda q H_b$ and hence $H_a=q H_b$.  But $q H_b=H_{q^{-1}b}$ and the equality $H_a=H_{q^{-1}b}$ together with Proposition \ref{propsubgroup}, $(i)$, show that $a=q^{-1}b$ in $ \A^f/\hatz$. Thus we get $\lambda'=q\lambda$ and $b=qa$
which proves the injectivity.
\endproof

In terms of semifields,   the space we are considering here is that of sub-semifields of $\R_{\rm max}$ which are abstractly isomorphic to algebraic extensions of $\Z_{\rm max}$. The natural action of $\R_+^\times$ on this space through the Frobenius automorphisms $\fr_u$ (at the multiplicative notation level, $\fr_u$ is raising to the power $u$) is given by $H\mapsto uH$.  Thus, the action of $\R_+^\times$ on the pair $(a,\lambda)$ is given by the formula
$
\fr_u(\Phi(a,\lambda))=\Phi(a,u\lambda).
$

\begin{thm} \label{structure3} The points of the arithmetic site  $\arith$ over $\rmax$ form the quotient of the ad\`ele class space of $\Q$ by the action of $\hatz$. The action of the Frobenius
automorphisms $\fr_\lambda$ of $\rmax$ on these points corresponds to the action of the id\`ele class group on the above quotient of the ad\`ele class space of $\Q$.
\end{thm}
\proof 
Consider a point of the arithmetic site over $\rmax$:  $f: (\pt,\rmax)\to\arith $. This is given by a pair $(p,f^\#_p)$ as in Definition \ref{defnpt}. By applying Theorem \ref{thmstart},  to a point $p$ of the topos $\wnt$ is associated a rank one ordered group $H$  and, by  Theorem \ref{structure2}, the stalk of the structure sheaf $\zmax$ at $p$ is the semifield $K=H_{\rm max}$.
The map of stalks  is given by a (local) morphism of semifields $f^\#_p:K\to \rmax$. One considers the following two cases, depending upon the nature of the range of $f^\#_p$. \newline
$(\alpha)$~The range of $f^\#_p$ is  the semifield $\B\subset \rmax$. Then, since $f^\#_p$ sends invertible elements to $1\in \B$, this map is completely determined and thus the pair $(p,f^\#_p)$ is uniquely determined by the point $p$. We write $\iota_p$ for such a degenerate morphism on the stalks: it exists uniquely for any point of $\wnt$. Thus, the  map  $p\mapsto (p,\iota_p)$ gives a bijection between points of $\wnt$ and  points of $\arith$ over $\rmax$ which are defined over $\B\subset \rmax$.  These points are the fixed points for the action of  the Frobenius
automorphisms $\fr_\lambda$ of $\rmax$.\newline
$(\beta)$~The range of $f^\#_p$ is not contained in  $\B$. In this case $f^\#_p$ is an isomorphism of the semifield $K$ with a sub-semifield of $\rmax$. Indeed, when restricted to the multiplicative group $K^\times$, $f^\#_p$ defines a non-trivial morphism of ordered groups and is  injective. Moreover, the range $f^\#_p(K)\subset \rmax$ is entirely determined by the non-trivial rank one subgroup $H'\subset \R$ defined by 
\begin{equation}\label{log}
H'=\{ \log u\mid u\in f^\#_p(K), u\neq 0\}.
\end{equation}
The subgroup $H'$ determines uniquely, up to isomorphism, the pair given by the point $p$ and the map $f^\#_p$. More precisely, if two pairs $f=(p,f^\#_p)$ and $g=(q,g^\#_q)$ define the same rank one subgroup $H'\subset \R$, then there exists a unique isomorphism of points $\phi:p\stackrel{\sim}{\to} q$ such that the following diagram commutes 
\[
\xymatrix@C=5pt@R=25pt{
\cO_p \ar[rr]^{\phi^\#} \ar[rd]_{f^\#_p} &&
\cO_{q} \ar[ld]^{g^\#_q} &
\\
&\rmax}
\]

This shows that the points of the arithmetic site  $\arith$ over $\rmax$ form the union of two sets: the set of isomorphism classes of points of $\wnt$ and the set of the non-trivial rank one subgroups $H\subset \R$. By Proposition \ref{propsubgroup} $(ii)$, the first set  coincides with the quotient
\begin{equation}\label{bc1}
\Q_+^\times\backslash\A^f/\hatz.
\end{equation}
By Lemma \ref{subbis}, the set of non-trivial rank one subgroups $H\subset \R$ is the quotient
\begin{equation}\label{bc2}
\Q_+^\times\backslash(\A^f\times \R_+^\times)/(\hatz\times \{1\})
\end{equation}
 of $(\A^f\times \R_+^\times)/(\hatz\times\{1\})$ by the diagonal action of $\Q_+^\times$. The quotient of the ad\`ele class space of $\Q$ by the action of $\hatz$
$$
\Q^\times\backslash\A/(\hatz\times \{1\})=\Q^\times\backslash(\A^f\times \R)/(\hatz\times \{1\})
$$
is also the disjoint union of two pieces. The first one arises from the contribution of $\A^f\times \{0\}\subset \A^f\times \R$ and gives $\Q^\times\backslash\A^f/\hatz=\Q_+^\times\backslash\A^f/\hatz$, since the element $-1\in \Q^\times$ is already in $\hatz$. Thus, this piece coincides with \eqref{bc1}. The second piece arises from the contribution of $\A^f\times \R^\times\subset \A^f\times \R$ and gives, using multiplication by $-1\in \Q^\times$ to shift back to positive reals, the quotient $\Q_+^\times\backslash(\A^f\times \R_+^\times)/(\hatz\times \{1\})$ which coincides with \eqref{bc2}. Thus, the points of the arithmetic site  $\arith$ over $\rmax$ coincide with the quotient of the ad\`ele class space of $\Q$ by the action of $\hatz$. The action of the Frobenius automorphisms $\fr_\lambda\in \Aut(\rmax)$ is the identity on the degenerate morphisms $\iota_p$ since these morphisms have range in the semifield $\B\subset \rmax$ which is fixed by the  $\fr_\lambda\in \Aut(\rmax)$. This result is in agreement with the fact that the scaling action of $\R_+^\times$ on $\Q^\times\backslash(\A^f\times \R)/(\hatz\times \{1\})$ fixes the contribution of $\A^f\times \{0\}\subset \A^f\times \R$. Under the correspondence given by \eqref{log}, the action of the Frobenius automorphisms $\fr_\lambda$ on the points translates into the scaling action $H\mapsto \lambda H$ on the rank one subgroups $H\subset \R$. In turns, by Lemma \ref{subbis}, this scaling action corresponds to the  action of the id\`ele class group $\GL_1(\A)/\GL_1(\Q)=\hatz\times \R_+^\times$ on the  quotient of the ad\`ele class space of $\Q$ by $(\hatz\times \{1\})$. Finally, notice that with $(p,f^\#_p)$ as in case $\beta)$, one gets in the limit $\lambda\to 0$, that $(p,\fr_\lambda\circ f^\#_p)\to (p,\iota_p)$. More precisely,  for any $x\in \cO_p$ one has $\fr_\lambda\circ f^\#_p(x)\to \iota_p(x)$. This is clear for $x=0$ (the neutral element for the addition), while for $x\neq 0$ and $u=f^\#_p(x)\in \rmax$ the action of $\fr_\lambda$ replaces $\log u$ by $\lambda \log u$ which converges to $0=\log 1=\log \iota_p(x)$ when $\lambda\to 0$.\endproof

\begin{rem}\label{casenote1}{\rm  Theorem 2.6 of \cite{CCas}  is  the analogue of  Theorem \ref{structure3} when the pair $\arith$ is replaced  by $(\wnt,\zminp)$. The proofs are similar. One replaces $\rmax$ in  $(\pt,\rmax)$ by the maximal compact subring $\cO(\rmax)$, where the functor $\cO$ is  defined in Remark \ref{casenote}.
}\end{rem}

\section{Hasse-Weil formula for the Riemann zeta function}\label{secthw}

In this section, we combine  Theorem \ref{structure3} with our previous results as in  \cite{CC1, CC2} and derive a description of the complete Riemann zeta function  as the Hasse-Weil zeta function of the arithmetic site. 

 \subsection{The periodic orbits of the Frobenius on points over $\rmax$}\label{periodicsect}

We shall describe, at a qualitative level, the periodic orbits of the Frobenius action on points of the arithmetic site $\arith$ over $\rmax$. These points contribute through various terms in the counting distribution. To obtain a conceptual formula, one considers the action of the id\`ele class group $G=\GL_1(\A)/\GL_1(\Q)=\hatz\times \R_+^\times$ on the ad\`ele class space of $\Q$. One then restricts the attention to test functions $h$ on $G$ which are invariant under the maximal compact subgroup $\hatz$, \ie of the form $h(u)=g(|u|)$, where $g$ is a test function on $\R_+^\times$ and $u\mapsto |u|$ is the module, \ie the projection of $G=\hatz\times \R_+^\times$ on the second factor.  Notice that restricting to $(\hatz\times \{1\})$-invariant test functions is the same as considering the action of $\R_+^\times$ on the  quotient of the ad\`ele class space of $\Q$ by $(\hatz\times \{1\})$. We shall also assume that  the support of the function $g$ is contained in
$(1,\infty)$. \newline
The periodic orbits of the Frobenius flow which contribute to the Lefschetz formula on the ad\`ele class space of $\Q$ correspond, at a place $v$ of $\Q$, to the ad\`eles $a=(a_w)$ which vanish at $v$, \ie with $a_v=0$. This condition makes sense in the  quotient of the ad\`ele class space of $\Q$ by $\hatz\times \{1\}$ and does not depend on the choice of a lift of a given class as an ad\`ele. Next, we investigate the meaning of this condition in terms of the point of $\arith$ over $\rmax$ that it labels, accordingly to the place $v$ of $\Q$.

\subsubsection*{Archimedean place}

For $v=\infty$ the archimedean place, the above condition means that the corresponding point $\alpha$ is in the contribution of \eqref{bc1} \ie it is fixed by the Frobenius action $\fr_\lambda(\alpha)=\alpha$ $\forall\lambda$. Thus, the points of $\arith$ over $\rmax$ which correspond to $v=\infty$ are the points of $\arith$ defined over $\B$.

\subsubsection*{Finite places} For a finite place $v=p$, with $p$  a rational prime, the condition $a_p=0$ implies that the class $\alpha$ of $a$ in the  quotient of the ad\`ele class space of $\Q$ by $(\hatz\times \{1\})$, fulfills the equation $\fr_\lambda(\alpha)=\alpha$  $\forall\lambda \in p^\Z$. Indeed, the principal ad\`ele $p:=(p_v)$ associated to the rational number $p$ is a product of the form
$
p=(u\times 1)(1\times p_\infty)p_p
$,
where each term of the product is an id\`ele. The first term belongs to $\hatz\times \{1\}$, the second to $1_{\A^f}\times \R_+^\times$ and the third term has all its components equal to $1$ except at $p$ where it is equal to $p$. The multiplication by this last term fixes $a$ since $a_p=0$. The multiplication by $(1\times p_\infty)$ is the action of $\fr_p$. Finally, the multiplication by $(u\times 1)$ is the identity in the  quotient of the ad\`ele class space of $\Q$ by $\hatz\times \{1\}$. Moreover, since one operates in the ad\`ele class space, the multiplication by the principal id\`ele $p$ is the identity, and thus one has $\fr_p(\alpha)=\alpha$. Conversely however, the condition $\fr_p(\alpha)=\alpha$ is not sufficient to ensure that the component $a_p=0$ since it holds for all fixed points of the Frobenius action. This condition is in fact best understood in terms of the intermediate semifields $\F\subset K\subset \overline{\F}=\Q_{\rm max}$ which label the points of $\wnt$. In these terms it means that the Frobenius map $x\mapsto x^p$ is an automorphism of $K$. Notice that this condition depends only on the point of the topos $\wnt $ over which the point of $\arith$ over $\rmax$ sits. This invariance of $K$ under the Frobenius at $p$ implies of course that the range of the  morphism $f^\#_p$ to $\rmax$ is invariant under $\fr_p$. The condition $\fr_p(\alpha)=\alpha$ is sufficient for non-degenerate points but for degenerate points (those defined over $\B$) it does not suffice to ensure  $a_p=0$.

\subsection{Counting function for  the Frobenius action on points over $\rmax$}\label{periodicsect}

In order to count the number of fixed points of the Frobenius action on points of $\arith$ over $\rmax$, we let $
\urep_w\xi(x)=\xi(w^{-1}x)
$ be the scaling action of the id\`ele class group $G=\GL_1(\A)/\GL_1(\Q)$  on the complex valued functions on the ad\`ele class space $\A_\Q/\Q^\times$ and we use the trace formula 
 in the form
\begin{equation}\label{geomside}
\Tr_{\rm distr}\left(\int h(w)\urep(w)d^*w\right )=\sum_v\int_{\Q^\times_v}\,\frac{h(w^{-1})}{|1-w|}\,d^*w.
\end{equation}
We refer to \cite{Co-zeta, Meyer, CMbook} for a detailed treatment.
The subgroups $\Q^\times_v\subset G=\GL_1(\A)/\GL_1(\Q)$ arise as isotropy groups.
One can understand why the terms $\displaystyle \frac{h(w^{-1})}{|1-w|}$ occur in the trace formula by computing, formally as follows, the trace  of the scaling operator $T=\urep_{w^{-1}}$
$$T\xi(x)=\xi(w x)=\int k(x,y)\xi(y)dy\,
 $$
 given by the distribution kernel $k(x,y)=\delta(w x-y)$
 $$
\Tr_{\rm distr}(T)=\int k(x,x)\,dx=\int \delta(w x-x)\,dx=\frac{1}{|w-1|}\int \delta(z)\,dz=\frac{1}{|w-1|}\,.
$$

 We apply the trace formula  \eqref{geomside}  by taking 
 the function $h$ of the form $h(u)=g(|u|)$, where the support of the function $g$ is contained in
$(1,\infty)$. On the left hand side of \eqref{geomside} one first performs the integration in the kernel $\hatz$ of the module $G\to \R_+^\times$. At the geometric level, this corresponds to taking the quotient of the ad\`ele class space by the action of $\hatz$. We denote by $\urep_w$ the scaling action on this quotient. By construction, this action only depends upon $|w|\in \R_+^\times$.
In order to count the number of fixed points of the Frobenius
 we consider the distributional trace of an expression
of the form
$
\int  g_x(u) \urep_u d^*u
$
where $g_x(u)=u\delta_x(u)$. It  is characterized, as a distribution, by its evaluation on test functions $b(u)$. This gives
\begin{equation}\label{deltafun}
    \int b(u)g_x(u)d^*u=b(x), \   \      \int  g_x(u) \urep_u d^*u=\urep_x
\end{equation}
 thus we are simply considering
an intersection number.
We now look at the right hand side of \eqref{geomside}, \ie at the terms
\begin{equation}\label{weil2}
\int'_{\Q_v^\times} \frac{h(w^{-1}) }{   \vert 1-w   \vert} \, d^* w.
\end{equation}
Since $h(w)=g(|w|)$ and  the support of the function $g$ is contained in
$(1,\infty)$, one sees that the integral \eqref{weil2} can be restricted  in all cases to the unit
ball $\{w\,; \,|w|<1\}$ of the local field $\Q_v$. Again we study separately the two cases.

\subsubsection*{Archimedean place}
At the archimedean place, let $u=w^{-1}$ then one has 
$$
\frac 12\left( \frac{1}{1-\frac 1u}+\frac{1}{1+\frac 1u}\right)=\frac{u^2}{u^2-1}.
$$
The above equation is applied for $u>1$, in which case one can write equivalently
\begin{equation}\label{weil1}
\frac 12\left( \frac{1}{\vert 1-u^{-1}\vert}+\frac{1}{\vert 1+u^{-1}\vert}\right)=\frac{u^2}{u^2-1}.
\end{equation}
The term corresponding to \eqref{weil2} yields the distribution $\kappa(u)$
\begin{equation}\label{kappadu}
\int_1^\infty\kappa(u)f(u)d^*u=\int_1^\infty\frac{u^2f(u)-f(1)}{u^2-1}d^*u+cf(1)\,, \qquad c=\frac12(\log\pi+\gamma)
\end{equation}
where $\gamma=-\Gamma'(1)$ is the Euler constant. The distribution $\kappa(u)$ is positive on $(1,\infty)$ where, by construction, it is given by $\kappa(u)=\frac{u^2}{u^2-1}$.
\vspace{.02in}

\subsubsection*{Finite places}

At a finite place one has
 $\vert 1-w   \vert=1$, thus for each finite prime $p$ one has
\begin{equation}\label{weil3}
\int'_{\Q_p^\times} \frac{h(w^{-1}) }{   \vert 1-w   \vert} \, d^* w=
\sum_{m=1}^\infty \log p \,\,g(p^m).
\end{equation}

 \subsection{Analogue of the Hasse-Weil formula for the  Riemann zeta function}\label{hasseweilsect}

The results of \S \ref{periodicsect} yield a distribution $N(u)$ on $[1,\infty)$ which counts the number of fixed points of the Frobenius action on points of $\arith$ over $\rmax$. The next task is to derive from this distribution the complete Riemann zeta function in full analogy with the geometric construction in finite characteristic. In that case, \ie  working over a finite field $\F_q$, the zeta function $\zeta_N(s)$, associated to a counting function $N(q^r)$ which counts the number of points over the finite field extensions $\F_{q^r}$, is given by the Hasse-Weil formula
\begin{equation}\label{zetadefn1}
 \zeta_N(s)=Z(q,q^{-s}), \  \  Z(q,T) = \exp\left(\sum_{r\ge 1}N(q^r)\frac{T^r}{r}\right) .
\end{equation}
Inspired by  the pioneering work of  R. Steinberg and J. Tits, C. Soul\'e  associated  a zeta function to any sufficiently regular counting-type function $N(q)$, by considering the following limit
\begin{equation}\label{zetadefn}
\zeta_N(s):=\lim_{q\to 1}Z(q,q^{-s}) (q-1)^{N(1)}\qquad s\in\R.
\end{equation}
Here,  $Z(q,q^{-s})$ denotes the evaluation, at $T=q^{-s}$,  of the  Hasse-Weil expression
\begin{equation}\label{zetadefn1}
Z(q,T) = \exp\left(\sum_{r\ge 1}N(q^r)\frac{T^r}{r}\right).
\end{equation}
For the consistency of the formula \eqref{zetadefn}, one requires that the counting function $N(q)$ is defined for all real numbers
 $q \geq 1$ and not only for prime integers powers as for the counting function in \eqref{zetadefn1}.  For many simple examples of rational algebraic varieties, like the projective spaces, the function  $N(q)$ is known to extend unambiguously  to all real positive numbers. The associated zeta function $\zeta_N(s)$ is easy to compute and it produces the expected outcome. For a projective line, for example, one finds  $\zeta_N(s) = \frac{1}{s(s-1)}$.
To by-pass the difficulty inherent to the definition \eqref{zetadefn}, when $N(1)=-\infty$,  one works  with the logarithmic derivative
\begin{equation}\label{normael}
    \frac{\partial_s\zeta_N(s)}{\zeta_N(s)}=-\lim_{q\to 1} F(q,s)
\end{equation}
where
\begin{equation}\label{fqsdefn}
F(q,s)=-\partial_s \sum_{r\ge 1}N(q^r)\frac{q^{-rs}}{r}\,.
\end{equation}
Then one finds, under suitable regularity conditions on $N(u)$, that (\cf~Lemma 2.1 of \cite{CC1})

\begin{lem} \label{compute1}With the above notations and for $\Re\mathfrak e(s)$ large enough, one has
\begin{equation}\label{lim}
    \lim_{q\to 1} F(q,s) = \int_1^\infty N(u)u^{-s}d^*u\,,\ \ d^*u=du/u
\end{equation}
and
\begin{equation}\label{logzetabis}
    \frac{\partial_s\zeta_N(s)}{\zeta_N(s)}=-\int_1^\infty  N(u)\, u^{-s}d^*u\,.
\end{equation}
\end{lem}
Using this lemma we can now recall one of our key results (\cf \cite{CC1, CC2})
\begin{thm}\label{mainthm} The zeta function associated by \eqref{logzetabis} to the counting distribution of \S \ref{periodicsect} is the  {\em complete} Riemann zeta function  $\zeta_\Q(s)=\pi^{-s/2}\Gamma(s/2)\zeta(s)$.
\end{thm}

This theorem implies that we are now working in a similar geometric framework as in the case of the geometry of a complete curve in finite characteristic. Moreover, in \cite{CC1} we also showed that the counting distribution $N(u)$ can be expressed in the form
 $$
N(u)=u-\frac{d}{du}\left(\sum_{\rho\in Z}{\rm order}(\rho)\frac{u^{\rho+1}}{\rho+1}\right)+1
$$
where $Z$ is the set of non-trivial zeros of the Riemann zeta function and the derivative is taken in the sense of  distributions. The main open question is at this point the definition of a suitable Weil cohomology which would allow one to understand the above equality as a Lefschetz formula.

\section{Relation of $\Spec\Z$ with the arithmetic site}\label{sectspz0}

As shown in \S \ref{periodicsect}, the periodic orbits of the Frobenius action which are responsible for the contribution of a given prime $p$ to the counting function correspond to adeles $(a_v)$ whose component $a_p=0$. This fact gives a first hint for the determination of  a precise relation between  $\Spec\Z$ and the arithmetic site. In this section we construct a geometric morphism of topoi $\Theta: \Spec\Z\to \wnto$, where $\Spec\Z$ is viewed as a topos with the Zariski topology. In the second part of the section we give an interpretation of the pullback of the structure sheaf of the arithmetic site in terms of Cartier divisors. 

\subsection{The geometric morphism $\Theta$}\label{sectspz}

Recall that $\Spec\Z$ is the affine scheme given by the algebraic spectrum of the ring $\Z$. As a  topological space it  contains, besides the generic point,  one point $p$ for each rational prime and the non-empty open sets are the complements of finite sets. The generic point belongs to all non-empty open sets. We consider the topos $ {\mathcal Sh}(\Spec\Z)$ of sheaves of sets on $\Spec\Z$ and investigate the definition of a geometric morphisms of topoi
$$
\Theta: \Spec\Z\to \wnto.
$$
The general theory states that any such morphism is uniquely determined by a flat functor 
$$
\Theta^*: \nto \longrightarrow {\mathcal Sh}(\Spec\Z).
$$
 Thus we are looking for a sheaf $\cS$ of sets on $\Spec\Z$ corresponding to the unique object $\ast$ of the small category $\nto$ and for an action of the semigroup $\nto$ on $\cS$. In fact, a geometric morphism such as $\Theta$ determines a map  between the points of $\Spec\Z$ and the points of $\wnto$.  We first give the natural guess for the definition of such map on points related to the periodic orbits associated to primes as in \S \ref{periodicsect}. For each prime $p$,  let $\alpha_p\in\A^f$  be the finite adele  whose components are all equal to $1$ except at $p$ where the component vanishes. We let $H(p)$ be the corresponding subgroup of $\Q$ as in \eqref{param} \ie
$$
H(p):=\{q\in \Q\mid \alpha_p q\in \hat\Z\}.
$$
Thus $H(p)$ is the group of fractions whose denominator is a power of $p$. To this ordered group  corresponds a point of $\wnt$ and also of $\wnto$. By working with $\wnto$ we see that, as a set with an action of $\nto$, the stalk of $\cS$ at the point $p\in \Spec\Z$ is $\cS_p=H(p)^+$, \ie the set of non-negative elements of $H(p)$ on which  $\nto$ acts by multiplication. Moreover $\wnto$ admits a point corresponding to the subgroup $\{0\}\subset \Q$: we associate this point to the generic point of $\Spec\Z$. Thus the stalk of $\cS$ at the generic point is $\{0\}$. This suggests that when one considers a section $\xi\in\Gamma(U,\cS)$ of  $\cS$ on an open set $U\subset\Spec\Z$, its value at all primes but finitely many is the zero element of the stalk $\cS_p$. \newline
In agreement with the above construction we introduce the following

\begin{defn} The sheaf of sets $\cS$ on $\spz$ assigns, for each open set $U\subset\Spec\Z$, the set $\Gamma(U,\cS)=\{U\ni p\mapsto \xi_p\in \cS_p=H(p)^+\vert\xi_p\neq 0~ \text{for finitely many primes} ~p\in U\}$. The action of $\nto$ on sections is done pointwise \ie $(n\xi)_p:=n\xi_p$.
\end{defn}

The restriction maps on sections of $\cS$ are the obvious ones and they are $\nto$-equivariant by construction. To check that the above definition produces a sheaf, we notice that for an open cover $\{U_j\}$ of a non-empty open set $U\subset\spz$, each $U_j\neq\emptyset$ has a finite complement so that the finiteness condition is automatic for the sections locally defined using that cover. Moreover, two sections $\xi, \eta\in\Gamma(U,\cS)$ with the same value at $p$ are equal on an open set containing $p$ and thus the stalk of $\cS$ at $p$ is $H(p)$.\newline
We now formulate the filtering condition for the functor $F:\nto\longrightarrow {\mathcal Sh}(\Spec\Z)$ that we just constructed. By applying Definition VII.8.1 of \cite{MM}, the three filtering conditions involve the notion of an epimorphic family of morphisms in the topos ${\mathcal Sh}(\Spec\Z)$.
Next lemma provides an easy way to prove that a  family of morphisms $f_i:\mathcal A_i\to \mathcal A$ is epimorphic.

\begin{lem}
Let $X$ be a topological space and ${\mathcal Sh}(X)$ a topos of sheaves (of sets) on $X$. Then, a  family of morphisms $f_i: \mathcal A_i\to \mathcal A$ is epimorphic if for any non-empty open set $U\subset X$ and a section $\xi\in \Gamma(U, \mathcal A)$ there exists an open cover $\{U_\alpha\}$ of $U$ and sections $\xi_\alpha\in \Gamma(U_\alpha, \mathcal A_{i(\alpha)})$ such that 
$$
f_{i(\alpha)}\circ \xi_\alpha=\xi\vert_{U_\alpha}.
$$
\end{lem}
\proof By definition, a  family of morphisms $f_i: \mathcal A_i\to \mathcal A$ is epimorphic if and only if for any morphisms $f,g: \mathcal A\to \mathcal B$ the equality $f\circ f_i=g\circ f_i$ for all $i$ implies $f=g$. Here, the hypothesis implies that for a section $\xi\in \Gamma(U, \mathcal A)$ there exists an open cover $\{U_\alpha\}$ of $U$ and sections $\xi_\alpha\in \Gamma(U_\alpha, \mathcal A_{i(\alpha)})$ such that
$$
f\xi\vert_{U_\alpha}=f\circ f_{i(\alpha)}\circ \xi_\alpha=g\circ  f_{i(\alpha)}\circ \xi_\alpha=g\xi\vert_{U_\alpha}
$$
so that one derives $f\xi=g\xi$ as required.\endproof

Alternatively (and equivalently) one can use the formalism of generalized elements which allows one to think in  similar terms as for the topos of sets. In the case of  the topos $\spz$ it is enough to consider the sections of the form $\xi\in\Gamma(U,\cS)$ where $U$ is an arbitrary open set, \ie a subobject of the terminal object $1$. We now reformulate the filtering condition as in Lemma VII.8.4 of \cite{MM}. In fact, we take advantage of the fact that the small category $\cC=\nto$ has a single object $\ast$ and that its image under the flat functor $F$ is the object $F(\ast)=\cS$ of ${\mathcal Sh}(\Spec\Z)$. 
\begin{enumerate}
\item For any open set $U$ of $\spz$ there exists a covering $\{U_j\}$ and sections $\xi_j\in \Gamma(U_j,\cS)$.
\item For any open set $U$ of $\spz$ and sections   $c,d\in\Gamma(U, \cS)$, there exists  a covering $\{U_j\}$ of $U$ and for each $j$ arrows $u_j,v_j:\ast \to \ast$ in $\nto$ and a section $b_j\in \Gamma(U_j,\cS)$ such that 
$$
c\vert_{U_j}=F(u_j)b_j, \  \  d\vert_{U_j}=F(v_j)b_j.
$$
\item Given two  arrows $u,v:\ast \to \ast$ in $\nto$ and any section $c\in\Gamma(U,\cS)$ with $F(u)c=F(v)c$, there
exists a covering $\{U_j\}$ of $U$ and for each $j$ an arrow $w_j:\ast \to \ast$,  and a section $z_j\in\Gamma(U_j,\cS)$, such that 
$$
F(w_j)z_j=c\vert_{U_j}, \  \  u\circ w_j=v\circ w_j\in \Hom_{\nto}(\ast,\ast).
$$
\end{enumerate}
We  now state the main result of this subsection 

\begin{thm}\label{thmspz}
The  functor $F:\nto\longrightarrow {\mathcal Sh}(\Spec\Z)$ which associates to the object $\ast$ the sheaf $\cS$ and to the endomorphisms of $\ast$ the natural action of $\nto$ on $\cS$ is filtering and defines a geometric morphism $\Theta: \Spec\Z\to \wnto$. The image of the point $p$ of $\spz$ associated to a prime $p$ is the point of $\wnto$ associated to the subgroup $H(p)\subset \Q$.
\end{thm}
\proof We check the three filtering conditions. Let $U$ be a non-empty open set of $\spz$. Then the $0$-section (\ie the section whose value at each $p$ is $0$) is a generalized  element of  $\Gamma(U,\cS)$. This checks $(i)$. We verify $(ii)$. Let $U$  be a non-empty open set of $\spz$ and $c,d$ be sections of $\cS$ over $U$. Then there exists a finite set $E\subset U$ of primes such that both $c$ and $d$ vanish in the complement of $E$. The complement $V:=U\setminus E$ is a non-empty open set of $\spz$ and for each $p\in E$  let $U_p=V\cup \{p\}\subset U$. By construction, the collection $\{U_p\}_{p\in E}$ form an open cover of $U$. The restriction of the sections $c,d$ of $\cS$ to $U_p$ are entirely determined by their value at $p$ since they vanish at any other point of 
$U_p$. Moreover, if one is given an element $b\in \cS_p$ in the stalk of $\cS$ at $p$, one can extend it uniquely to a section of $\cS$ on $U_p$ which vanishes in the complement of $p$. Thus, for each $p$ one checks, using the filtering property of the functor associated to the stalk $\cS_p$, that, as required, there exist arrows $u_p,v_p\in\nto$ and a section $b_p\in \Gamma(U_p,\cS)$ such that 
$
c\vert_{U_p}=F(u_p)b_p, \   d\vert_{U_p}=F(v_p)b_p.
$
Since $\{U_p\}_{p\in E}$ form an open cover of $U$, one  gets $(ii)$. Finally, we verify $(iii)$. Let $c\in\Gamma(U,\cS)$ with $F(u)c=F(v)c$. Assume first that for some $p\in U$ one has $c_p\neq 0$. Then the equality $F(u)c=F(v)c$ in the stalk at $p$ implies that $u=v$ and in that case one can take the cover of $U$ by $U$ itself and $z=c$, $w=1$. Otherwise,  $c$ is the zero section and one can take the cover of $U$ made by $U$ itself and $z=0$, $w=0$. We have thus shown that the functor $F:\nto\longrightarrow {\mathcal Sh}(\Spec\Z)$ is filtering and by Theorem VII.9.1 of \cite{MM} it is thus flat so that by Theorem VII.7.2 of \opcit 
it defines a geometric morphism $\Theta: \Spec\Z\to \wnto$. The image of a point $p$ of $\spz$ is the point of $\wnto$ whose associated flat functor $G:\nto\longrightarrow \Se$ is obtained by composing the functor $F:\nto\longrightarrow {\mathcal Sh}(\Spec\Z)$
with the stalk functor at $p$. This functor assigns to any sheaf on $\spz$ its stalk at $p$, viewed as a set. In the above case we see that $G$ is the flat functor $G:\nto\longrightarrow \Se$ given by the stalk $\cS_p=H(p)^+$ and this proves the  statement.\endproof

\begin{prop}\label{sheavesonspz0} Let $\Theta^*(\cO)$ denote the pullback of the structure sheaf of $\arith$. Then

$(i)$~The stalk of $\Theta^*(\cO)$ at the prime $p$ is the semiring $H(p)_{\rm max}$, and it is 
$\B$ at the generic point of $\Spec\Z$.

$(ii)$~The sections $\xi\in\Gamma(U,\Theta^*(\cO))$ on an open set $U\subset\spz$ are the maps $U\ni p\mapsto \xi_p\in H(p)\subset H(p)_{\rm max}$ which are either equal to $0$ outside a finite set,  or everywhere equal to the constant section   $\xi_p=-\infty \in H(p)_{\rm max}$, $\forall p\in U$.
\end{prop}

\proof $(i)$~The stalk of $\Theta^*(\cO)$ at the prime $p$ is the same as the stalk of $\cO$ at the point $\Theta(p)$ of $\wnto$. Thus $(i)$ follows from Remark \ref{caseno}.\newline
$(ii)$~The pullback functor $ \Theta^*$ commutes with arbitrary colimits. As a sheaf of sets on $\wnto$ the structure sheaf $\cO$ is the coproduct of two sheaves. The first is the constant sheaf taking always the value $-\infty$. The second is the coproduct of two copies of $\nto$ glued over the map from the constant sheaf $0$ to the zeros in both terms. It  follows that the sheaf $ \Theta^*(\cO)$ on $\spz$ is obtained as the corresponding colimit among sheaves on $\spz$. One has  $ \Theta^*(\nto)=\cS$ and the pullback of the constant sheaves  $-\infty$ and $0$ are the corresponding constant sheaves on $\spz$. Thus $ \Theta^*(\cO)$ is the coproduct of the constant sheaf $-\infty$ on $\spz$ with the sheaf on $\spz$ obtained by glueing two copies $\cS_\pm$ of $\cS$ over the constant sheaf $0$ which maps to $0_\pm$. Let $U$ be a non empty open set of $\spz$ and $U\ni p\mapsto \xi_p\in H(p)\subset H(p)_{\rm max}$ a map equal to $0$ outside a finite set $E\subset U$. To show that $\xi\in\Gamma(U,\Theta^*(\cO))$, let $E_+=\{p\mid \xi_p>0\}$ and $E_-=\{p\mid \xi_p<0\}$. Consider the sections $\xi_\pm=\xi\vert_{U_\pm}\in\Gamma(U_\pm,\cS_\pm)$ obtained as restrictions of $\xi$ to $U_+=U\setminus E_-$ and $U_-=U\setminus E_+$. One has $\xi_+=\xi_-$ on $U_+\cap U_-$ since one identifies the $0_+$ section with the $0_-$ section on any open set. Thus since $ \Theta^*(\cO)$ is a sheaf, there is a section $\eta \in\Gamma(U,\Theta^*(\cO))$ which agrees with $\xi_\pm$ on $U_\pm$ and hence $\eta=\xi$ and $\xi\in\Gamma(U,\Theta^*(\cO))$. Finally there is only one section of the constant sheaf $-\infty$ on any non-empty open set of $\spz$. 
\endproof 

\subsection{$\Theta^*(\cO)$ as Cartier divisors on $\Spec\ocy$}
In this part, we relate the sheaf $\Theta^*(\cO)$ with Cartier divisors on the algebraic spectrum  of the ring $\ocy$ of the cyclotomic integers. By construction, $\ocy$ is a filtering limit of the rings of integers $\cO_K$ of finite, intermediate field extensions $\Q\subset K\subset \qcy$. It is a well-known fact that the localization of $\cO_K$ at any non-zero prime ideal is a discrete valuation ring (\cf \cite{Hart} \S I.6). 
To the ring homorphism  $\Z\to  \cO_K$ corresponds a morphism of schemes $\pi_K:\Spec\cO_K\to \Spec\Z$ and for a prime ideal $\wp\subset \cO_K$ above a given rational prime $p$, one  normalizes the associated discrete valuation  so that $v_\wp(p)=p-1$.\newline
We recall that the Cartier divisors on a scheme $X$ are the global sections of the sheaf $\cK^\times/\cO^\times$ which is the quotient of the sheaf of multiplicative groups of the total quotient rings $\cK$ of $X$ by the sub-sheaf $\cO^\times$ of the structure sheaf (\cf \cite{Hart} page 140-141). Since $\cK$  is a sheaf of rings and $\cO^\times$ is a sheaf of subgroups of the multiplicative groups one gets on the quotient $\cK/\cO^\times$ a natural structure of sheaf of hyperrings.

\begin{defn}\label{cartierdiv} Let $X$ be an algebraic scheme. We denote by $\card(X)$ the sheaf of hyperrings $\cK/\cO^\times$.
\end{defn}

Since $X=\Spec\cO_K$ is an integral scheme, $\cK$ is simply a constant sheaf $K$. 

\begin{lem}\label{hyper}
Let $\wp\subset\cO_K$ be a prime ideal above a rational prime $p\subset\Z$. Then, the normalized valuation $v_\wp$ determines an isomorphism of hyperfields 
\begin{equation}\label{minusvj}
-v_\wp:K/(\cO_K)_\wp^\times\stackrel{\sim}{\to} \alpha_\wp\zhmax,\quad \alpha_\wp\in\Q_{>0}
\end{equation}
where $\zhmax$ is the hyperfield obtained from the semifield $\zmax$ by altering the addition of two equal elements $x\in\zmax$ as
\begin{equation}\label{alterh}
x\vee x:=[-\infty,x].
\end{equation}
\end{lem}
\proof By definition one has $v_\wp(0)=\infty$. For $x,y\in K$ one has: $x\in (\cO_K)_\wp^\times y \iff v_\wp(x)=v_\wp(y)$. Thus $-v_\wp$ is a bijection of   $K/(\cO_K)_\wp^\times$ with  the range $\alpha_\wp\Z\cup \{-\infty\}$ of $-v_\wp$, where $\alpha_\wp$ is a positive rational number. This bijection transforms the multiplication of $K$ in addition. Let then $x,y\in K$ and assume that $v_\wp(x)\neq v_\wp(y)$, then $v_\wp(x+y)=\inf (v_\wp(x),v_\wp(y))$ so that one gets the additive rule 
$$
-v_\wp(x+y)=-v_\wp(x)\vee -v_\wp(y).
$$
For $v_\wp(x)= v_\wp(y)=-a$ one checks that all values in $[-\infty,a]\subset \alpha_\wp \zmax$ are obtained as $-v_\wp(x+y)$. Since  the quotient of a ring by a subgroup of its multiplicative group defines always a hyperring one knows that $\zhmax$ is an hyperfield. \endproof

\begin{prop}\label{sheavesonspz} The normalized valuation $-v_\wp$ determines an isomorphism of sheaves of hyperrings on $\Spec(\ocy)$
\begin{equation}\label{relsheaf2}
\card(\Spec(\ocy))\stackrel{\sim}{\to} (\Theta\circ \pi)^*(\cO_h)
\end{equation}
where $\cO_h$ denotes the sheaf of hyperrings on the arithmetic site obtained from the structure sheaf by altering the addition as in \eqref{alterh}.
\end{prop}
\proof  At the generic point $\xi$ of $\Spec\cO_K$ one has $(\cO_K)_\xi=K$ and the stalk of $\card(\Spec\cO_K)$ is the Krasner hyperfield $\K=\{0,1\}$ ($1+1=\{0,1\}$). This coincides with  $(\Theta\circ \pi)^*(\cO_h)$ since the stalk of the structure sheaf $\cO$ of the arithmetic site at the point $\Theta(\pi(\xi))$ is the semiring $\B$. Thus at the generic point,  \eqref{relsheaf2} is an isomorphism. The stalk of $(\Theta\circ \pi)^*(\cO_h)$ at a point $\wp\neq \xi$ is the hyperfield  $H(p)_{\rm hmax}$ obtained from the semiring $H(p)_{\rm max}$ by altering the addition as in \eqref{alterh}.

Let $v\in{\rm Val}_p(\qcy)$ be an extension of the $p$-adic valuation to $\qcy$. Let us check that with the normalization taken above, \ie $v(p)=p-1$, one has $v(\qcy)=H(p)$.  The composite subfield $\Q_p \vee \qcy\subset(\qcy)_v$  is the maximal abelian extension $\Q_p^{\rm ab}$ of $\Q_p$.  This extension  coincides  with the composite 
$
    \Q_p^{\rm ab}=\Q_p^{\rm ur}\vee \Q(\mu_{p^\infty})
$, 
where $\Q_p^{\rm ur}$ denotes the maximal unramified extension of $\Q_p$ and $\Q(\mu_{p^\infty})$ is obtained by adjoining to $\Q_p$ all roots of unity of order a $p$-power (\cf \cite{Serre}). The extension $\Q_p^{\rm ur}$ does not increase the range of the valuation since it is unramified.  The extension $ \Q(\mu_{p^\infty})$ is totally ramified and  the $p$-adic valuation extends uniquely from $\Q$ to $ \Q(\mu_{p^\infty})$.  Let $z$ be a primitive root of $1$ of order $p^m$. Then $u=z^{p^{m-1}}$ fulfills the equation $u^{p-1}+u^{p-2}+\cdots+ u+1=0$, and with  $z=1+\pi$  the equation fulfilled by $\pi$ is of Eisenstein type, the constant term being equal to $p$,  and reduces to $\pi^{\varphi(n)}=0$,  modulo $p$. This shows that
$
v(\pi)=\frac{v(p)}{\varphi(n)}\,, \ \varphi(n)=(p-1)p^{m-1}.
$
Thus with the valuation $v$, normalized as $v(p)=p-1$, 
 one gets $v(\qcy)=H(p)$ as required. \newline
Let $\Q\subset K\subset \qcy$ be a finite extension of $\Q$ and $\cO_K$ its ring of integers. The non-empty open sets $U\subset \Spec\cO_K$ are the complements of finite sets of prime ideals. We claim that the map which to $u\in K$ associates $\phi(u)=\xi$, with $\xi_\wp:=-v_\wp(u)\in H(\pi(\wp))_{\rm max}$, for any prime ideal $\wp$ of  $\cO_K$, defines a morphism of sheaves $\card(\Spec(\cO_K))\to (\Theta\circ \pi)^*(\cO_h)$. For $u=0$, the image of $\phi$ is the constant section $-\infty$ which is  the pullback by $\pi$ of the constant section $-\infty$ on $\spz$ as in Proposition \ref{sheavesonspz0}. For $u\neq 0$, one needs to show that  the section $\phi(u)=\xi$ is a section of the pullback by $\pi$ of the sheaf $\Theta^*(\cO)$ again as in  Proposition \ref{sheavesonspz0}. To this end, it suffices to find an open cover $\{U_j\}$ of $\Spec\cO_K$
and for each $U_j$ a section $\eta_j\in\Gamma(\spz,\Theta^*(\cO)) $ such that $\xi\vert_{U_j}=
(\pi^*\eta_j)\vert_{U_j}$. Since $u\neq 0$, one has $u\in K^\times$ and thus  $-v_\wp(u)=0$ outside a finite set $F$ of prime ideals of $\cO_K$. Let $F'=\pi^{-1}\pi(F)$ be the  saturation of $F$: it is a finite subset of  $\Spec\cO_K$ and its complement $U$ is open. For each $j\in F'$ one lets $U_j=U\cup\{j\}$. One has by construction $-v_\wp(u)\in H(\pi(j))$ and  by Proposition \ref{sheavesonspz0}, there is a  section $\eta_j$ of the sheaf $\Theta^*(\cO)$ which takes the value $-v_j(u)$ at $\pi(j)$ and vanishes elsewhere. One then gets $\xi\vert_{U_j}=
(\pi^*\eta_j)\vert_{U_j}$ as required. This reasoning applied in junction with Lemma \ref{hyper} show  that the map $\phi$ defines an injective morphism of sheaves $\card(\Spec(\cO_K))\to (\Theta\circ \pi)^*(\cO_h)$. These morphisms are compatible with the filtering colimit  expressing $\ocy$ in terms of subrings $\cO_K$ and in the colimit one obtains the surjectivity at the level of the stalks and hence the required isomorphism. \endproof

\section{The square of the arithmetic site and  Frobenius correspondences}\label{sectsquare}

The last sections have provided an interpretation of the non-commutative space quotient of the ad\` ele class space of $\Q$ by the action of $\hatz$ as the set of points of the arithmetic site over  $\rmax$. The action of $\R_+^\times$ on the ad\`ele class space corresponds to the action of the Frobenius $\fr_\lambda\in \Aut(\rmax)$ on the above points. It is natural to ask what one gains with this interpretation with respect to the original non-commutative geometric approach. In this section we prove that this algebro-geometric construction  leads naturally to the definition of the square of the arithmetic site and to the introduction of the Frobenius correspondences.

\subsection{The closed symmetric monoidal category $\catmo(\B)$}

We refer to \cite{PR} for a general treatment of tensor products of semirings. The category $\catmo(\B)$ of $\B$-modules endowed with the tensor product  is a closed, symmetric, monoidal category.  One first defines the internal homomorphisms  
$\underline\Hom$ by endowing the set of homomorphisms of two $\B$-modules $\Hom_\B(E,F)$ with the following addition
$$
(\phi+\psi)(x):=\phi(x)+\psi(x)\qqq x\in E.
$$
The definition of the tensor product is then dictated by the required adjunction formula
$$
\Hom_\B(E\otimes_\B F,G)\simeq \Hom_\B(E,\underline\Hom(F,G)).
$$
The tensor product $E\otimes_\B F$ of $\B$-modules is, as in the case of modules over rings,  the (unique) initial object in the category of bilinear maps $E\times F\to G$. It is constructed as  the quotient  of the $\B$-module of finite formal sums $\sum e_i\otimes f_i$ (where no coefficients are needed since $E$ and $F$ are idempotent) by the equivalence relation 
$$
\sum e_i\otimes f_i\sim \sum e'_j\otimes f'_j\iff \sum \rho(e_i,f_i)=\sum \rho(e'_i,f'_i)\qqq \rho
$$
where $\rho$ varies through all bilinear maps from $E\times F$ to arbitrary $\B$-modules. 
\begin{rem}  {\rm The notion of tensor product reviewed above does not coincide with the tensor product considered in Chapter 16 of \cite{Golan}. Indeed, it follows from Proposition (16.12) of \opcit that any tensor product so defined is cancellative which excludes all $\B$-semimodules $M$ since for any $x,y\in M$ one has $x+z=y+z$ for $z=x+y$. The notion recalled above is developed \eg  in \S 6 of \cite{PR}, and had been previously well understood by \eg G. Duchamp (letter to P. Lescot).
}\end{rem}

\subsection{The semiring $\nbo$ and the unreduced square $\arithb$}

In this subsection we describe the  tensor product 
$
\zmin\otimes_\B \zmin
$.
Let $R$ be a $\B$-module and consider bilinear maps $\phi: \zmin\times \zmin\to R$, \ie maps  fulfilling
\begin{equation*}
\phi(a\wedge b,c)=\phi(a,c)\oplus\phi(b,c), \  \  \phi(a, b\wedge c)=\phi(a,b)\oplus\phi(a,c), \  \  \phi(\infty,b)=\phi(a,\infty)=0
\end{equation*}
where  $\oplus$ denotes the idempotent addition of $R$ and $0\in R$ the neutral element. The simplest example of such a map  is the projection $\pi: \zmin\times \zmin\to\zmin\wedge\zmin$, where the smash product is taken for the base point $\infty\in \zmin$ and the idempotent addition is 
defined component-wise as $(a,b)\oplus(c,d):=(a\wedge c, b\wedge d)$ outside the base point, \ie for $a,b,c,d\in\Z$. The bilinearity follows using the idempotent property of addition in $\zmin$. At first, one might think  that any bilinear map $\phi: \zmin\times \zmin\to R$ would factor through this $\B$-module $\zmin\wedge \zmin$, however this is not true even for the operation $m(a,b)= a+b$ in $\zmin$ playing the role of the product in the semiring structure. One might then conclude from this example that one would need to consider all finite formal sums of simple tensors $a \otimes_\B b$, for $a,b\in \zmin$, \ie roughly, the finite subsets of $\Z^2$. However, the associativity of the idempotent addition $\oplus$ of $R$ provides one further relation as the following lemma points out

\begin{lem} \label{assoc} For any bilinear map $\phi: \zmin\times \zmin\to R$  one has 
$$
x_1\wedge x_2=x_1, \  y_1\wedge y_2=y_1\implies \phi(x_1\wedge x_2,y_1\wedge y_2)=\phi(x_1,y_1)\oplus\phi(x_2,y_2).
$$
\end{lem}
\proof Using the bilinearity of $\phi$ and the hypothesis one has
$$
\phi(x_1,y_1)\oplus\phi(x_2,y_2)=\phi(x_1\wedge x_2, y_1)\oplus \phi(x_2,y_2)=\left(\phi(x_1,y_1)\oplus\phi(x_2,y_1)\right)\oplus \phi(x_2,y_2)=
$$
$$
=\phi(x_1,y_1)\oplus\left(\phi(x_2,y_1)\oplus \phi(x_2,y_2)\right)=
\phi(x_1,y_1)\oplus \phi(x_2,y_1\wedge y_2)=
=
\phi(x_1,y_1)\oplus \phi(x_2,y_1)$$
$$=
\phi(x_1\wedge x_2,y_1)=\phi(x_1\wedge x_2,y_1\wedge y_2).
$$ \endproof
The above result  shows that the following rule  holds in $\zmin\otimes_\B \zmin$
\begin{equation}
x_1\wedge x_2=x_1, \  y_1\wedge y_2=y_1\implies (x_1\otimes_\B y_1)\oplus (x_2\otimes_\B y_2)=x_1\otimes_\B y_1.
\end{equation}

In view of these facts we introduce the following

\begin{defn} Let $J$ be a partially ordered set. We denote by $\sub(J)$ the set of subsets $E\subset J$ which are finite unions of intervals of the form $I_x:=\{y\mid y\geq x\}$.
\end{defn}

Such subsets $E\subset J$ are  hereditary, \ie they fulfill the rule:  $x\in E\implies y\in E$, $\forall y\geq x$.

\begin{lem} \label{assoc1} Let $J$ be a partially ordered set, then $\sub(J)$ endowed with the operation $E \oplus E':=E\cup E'$ is a $\B$-module.
\end{lem}
\proof One just needs to check that $E\cup E'\in \sub(J)$ which is immediate. The empty set $\emptyset \in \sub(J)$ is the neutral element for the operation $\oplus$. \endproof

\begin{lem} \label{assoc2} Endow $J=\N\times \N$ with the product ordering:  
$(a,b)\leq (c,d)\iff a\leq c \  \text{and} \  b\leq d$.
Then, $\sub(\N\times \N)$ is the set of hereditary subsets of $\N\times \N$. \end{lem}
\proof Let $E\neq \emptyset$ be a hereditary subset of $\N\times \N$ and let $(a,b)\in E$. Then, for any $x\in \N$ let   $L_x:=\{y\in \N\mid (x,y)\in E\}$. One has $L_x\subset L_z$ for $x\leq z$ (because $(x,y)\in E$ implies $(z,y)\in E$). Thus there exists $x_{\rm min}\in [0,a-1]$ such that 
$
L_x\neq \emptyset \iff x\geq x_{\rm min}
$.
Let  $\gamma: [x_{\rm min},\infty)\to \N$ be such that $L_x=[\gamma(x),\infty)$ for all $x\geq x_{\rm min}$. The map $\gamma$ is non-increasing and thus $\gamma$ takes only finitely many values. We thus get finitely many pairs $\alpha_j=(x_j,y_j)$ such that $\gamma(x_j)<\gamma(x_j-1)$. By construction of $\gamma$, the set $E$ is the upper graph of $\gamma$, and one gets $E\in \sub(\N\times \N)$ since $E$ is the finite union of the intervals  $I_{(x_j,y_j)}$.
\endproof
\begin{prop}\label{tensjustif} Let $S=\sub(\Z\times \Z)$ be the $\B$-module associated by Lemma \ref{assoc1} to the ordered set $\Z\times \Z$ with the product ordering:  
$(a,b)\leq (c,d)\iff a\leq c \  \text{and} \  b\leq d$. Then\newline
$(i)$~The following equality defines a bilinear map: 
\begin{equation}\label{eset1}
\psi:\zmin\times \zmin\to S,\quad \psi(u,v)=\{(a,b)\in \Z\times \Z\mid a\geq u, \, b\geq v\}.
\end{equation}

$(ii)$~Let $R$ be a $\B$-module and  $\phi: \zmin\times \zmin\to R$ be bilinear. Then there exists a unique $\B$-linear map $\rho: S\to R$ such that $\phi=\rho\circ \psi$. In other words
one has the identification 
\begin{equation}\label{square1}
\zmin\otimes_\B \zmin=\sub(\Z\times \Z).
\end{equation}
\end{prop}
\proof $(i)$~For a fixed $v\in \Z$, the map $u\mapsto \psi(u,v)$ which associates to $u\in \Z$ the quadrant as in \eqref{eset1} is monotone. This fact gives the required linearity.\newline
$(ii)$~By construction, an element of $S$ is a finite sum of the form $z=\oplus \psi(\alpha_i)$. The uniqueness of $\rho$ is implied by  $\phi=\rho\circ \psi$ which shows that
$
\rho(z)=\sum \phi(x_i,y_i), ~\alpha_i=(x_i,y_i)~ \forall i
$.
To prove the existence of $\rho$, we show  that the sum $\sum \phi(x_i,y_i)$ does not depend on the choice of the decomposition $z=\oplus \psi(\alpha_i)$. In fact in such a decomposition the set of the $\alpha_i$'s necessarily contains as a subset the elements $\beta_\ell$ of the canonical decomposition of Lemma \ref{assoc2}. Moreover, for any element $\alpha_j$ which is not a $\beta_\ell$, there exists a $\beta_\ell$ such that $\alpha_j\geq \beta_\ell$ for the ordering of $\Z \times \Z$. Thus one concludes, using Lemma \ref{assoc}, that one has: 
$
\sum \phi(x_i,y_i)=\sum \phi(a_\ell,b_\ell),   \   \beta_\ell=(a_\ell,b_\ell)
$
independently of the decomposition $z=\oplus \psi(\alpha_i)$. It then follows that $\rho$ is additive and fulfills $\phi=\rho\circ \psi$. \endproof

\begin{figure}
\begin{center}
\includegraphics[scale=0.4]{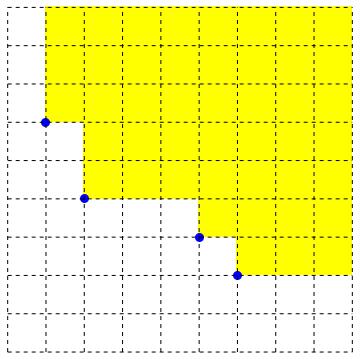}
\caption{Typical subset $E\in \sub(\N\times \N)$}\label{adsub}
\end{center}
\end{figure}

From now on, we shall  adopt  a multiplicative notation
 to refer to elements in $\zmin$, \ie we associate to $n\in\zmin$ the {\em symbol} $q^n$. The second operation of $\zmin$ then becomes the ordinary product 
and if one represents $q$ as a real number $0<q<1$, the first operation corresponds to the addition in $\rmax$: $x\vee y:=\max(x,y)$.  Here, $\rmax$ denotes the multiplicative version of the tropical semifield $\R_{\rm max}:= (\R\cup\{-\infty\},\text{max},+)$.

Every element of $\zmin\otimes_\B \zmin$ is a finite sum of the form
$
x=\sum q^{n_i}\otimes_\B q^{m_i}
$.
 We can endow the tensor product $S=\zmin\otimes_\B \zmin$ with the structure of a semiring of characteristic $1$. 
\begin{prop}\label{tensring} 
$(i)$~On the $\B$-module  $S=\zmin\otimes_\B \zmin$ there exists a unique bilinear multiplication satisfying the rule 
\begin{equation}\label{square3}
(q^ a\otimes_\B q^b)(q^ c\otimes_\B q^d)=q^ {a+c}\otimes_\B q^{b+d}.
\end{equation} 
$(ii)$~The above multiplication turns $S=\zmin\otimes_\B \zmin$ into a semiring of characteristic $1$.\newline
$(iii)$~The following equality defines an action of $\N^\times\times \N^\times$  by endomorphisms of $\zmin\otimes_\B \zmin$
\begin{equation}\label{bifrob}
\fr_{n,m}(\sum q^a\otimes_\B q^b):=\sum q^{na}\otimes_\B q^{mb}.
\end{equation}
\end{prop}
\proof We define an operation directly on $\sub(\Z\times \Z)$ as follows
\begin{equation}\label{square2}
E +E':=\{\alpha+\alpha' \mid \alpha \in E, \  \alpha' \in E'\}.
\end{equation}
One has $\psi(u,v)=(u,v)+(\N\times \N)$ and thus 
$
\psi(u,v)+\psi(u',v')=\psi(u+u',v+v')
$. 
Since $\psi(u,v)$ corresponds to $q^ u\otimes_\B q^v$ under the identification \eqref{square1}, one gets $(i)$. This also shows that the operation \eqref{square2} is well defined, \ie that $E +E'\in \sub(\Z\times \Z)$.
Moreover $\psi(0,0)=(\N \times \N)$ is the neutral element. The operation \eqref{square2} is associative and commutative and one has $(E \cup E')+E''=(E+E'')\cup(E'+E'')$ so that one obtains a semiring structure on 
$\sub(\Z\times \Z)$ and one gets $(ii)$.\newline
$(iii)$~The affine transformation of $\R^2$ given by the matrix $\left(
\begin{array}{cc}
 n & 0 \\
 0 & m \\
\end{array}
\right)$ maps $\N\times\N$ to itself but fails to be surjective. It is an automorphism of the quadrant $Q:=\R_+\times \R_+$ (in fact also of $\Q_+\times \Q_+$). Let us associate to any $E\in \sub(\Z\times\Z)$ the subset of  $\R^2$ given by $\tilde E:=E+Q$. One then has $E=\tilde E\cap (\Z\times \Z)$  so that the map $E\mapsto\tilde E$ is injective. Moreover the two operations in $\sub(\Z\times\Z)$ satisfy
$
\widetilde{ (E\cup E')}=\tilde E \cup \tilde E', ~ \widetilde{ (E+ E')}=\tilde E + \tilde E'
$.
 The action \eqref{bifrob} is then given in terms of the associated subsets by 
$
\widetilde{\fr_{n,m}(E)}=\left(
\begin{array}{cc}
 n & 0 \\
 0 & m \\
\end{array}
\right)\, \tilde E
$
which is simply an affine transformation preserving $Q$ and commuting with the operations of union and sum.
 \endproof

\begin{rem}\label{Steph}{\rm We are grateful to S. Gaubert, who pointed out to us that 
 the semiring $\zmin\otimes_\B\zmin$
has been introduced in a totally different context, namely in the modelisation 
of execution of discrete events.
In that framework, it is denoted  $\cM_{\rm in}^{\rm ax}[[\gamma,\delta]]$ and corresponds more precisely to $\zmin\otimes_\B\zmax$.
In that set-up  $q^a \otimes_B q^b$ encodes the fact that the event labelled  $a$ 
occurs  after the instant labelled $b$. We refer to \cite{Gaubert1} for the convex analysis and spectral analysis of timed event graphs.
}\end{rem}

Proposition \ref{tensring} shows that the semiring $\nbo$ is enriched with an action of the multiplicative mono\"id $\nt\times \nt$  of pairs of non-zero positive integers. This action is given by the endomorphisms $\fr_{n,m}\in \End(\nbo)$
\begin{equation}\label{Fr}
\nt\times \nt \to \End(\nbo),\qquad (n,m)\to \fr_{n,m}
\end{equation}
\begin{rem}{\rm Notice that the diagonal $\fr_{n,n}\in \End(\nbo)$ does not coincide with the operation 
$x\mapsto x^n$ in $\nbo$. In fact this operation fails to be an endomorphism of $\nbo$ which itself fails to be multiplicatively cancellative. This defect will be taken care in the following part of this section by passing to the reduced semiring.
}\end{rem}

We denote by $\wntb$ the topos of sets endowed with an action of  $\nt\times \nt$. It is the dual of the small category with a single object $\{*\}$ and whose endomorphisms form the semigroup $\nt\times \nt$.

\begin{defn}\label{site2} The {\em unreduced square} $\arithb$ of the arithmetic site is the topos $\wntb$ with the  {\em structure sheaf} $\nbo$, viewed as a semiring in the topos.
\end{defn}

The unreduced square comes endowed with the two projections $\pi_j$ to $\arith\sim (\wnt,\zmin)$. These morphisms are specified by the corresponding geometric morphisms of topoi and the semiring homomorphisms. We describe the projection $\pi_1$ on the first factor. One defines $p_1:\N^{\times 2}\to \N^\times$ as $p_1(n,m):=n$: it gives a corresponding geometric morphism $\hat p_1: \wntb\to \wnt$. At the semiring level, one has a natural homomorphism $\iota_1:\zmin\to \nbo$ given in multiplicative notation by $q^n\mapsto q^n\otimes 1$. Moreover, one derives the compatibility:
$
\fr_{n,m}(\iota_1(x))=\iota_1(\fr_{p_1(n,m)}(x))\ \forall x\in \zmin.
$

\begin{prop}\label{tensjustifsquare} The set of points of the unreduced square $\arithb$ over $\rmax$ is canonically isomorphic to the square of the set of points of $\arith$ over $\rmax$.
\end{prop}
\proof By Definition \ref{defnpt}, a point of the unreduced square $\arithb$ over $\rmax$ is given by a pair of a point $p$ of $\wntb$ and a  morphism $f^\#_p:\cO_p\to \rmax$ of the stalk of the structure sheaf at $p$ to $\rmax$. As a multiplicative mono\"{\i}d, the product $\nt\times \nt$ is isomorphic to $\nt$ and it follows that, up to isomorphism, all its points are obtained as pairs of the points of $\wnt$ associated to rank one ordered groups. Moreover the stalk of the structure sheaf at the point $p$ of $\wntb$ associated to the pair $(H_1,H_2)$ of rank one ordered groups is given by the tensor product of semirings $\cO_p=H_{1,{\rm max}}\otimes_\B H_{2,{\rm max}}$. Finally a morphism of semirings $f^\#_p:\cO_p\to \rmax$ is the same thing as a pair of morphisms $H_{j,{\rm max}}\to \rmax$.\endproof

\subsection{The Frobenius correspondences on $\arith$}\label{sectfrob}

Next, we want to interpret the product in the semiring $\zmin$ as yielding a morphism of semirings $\mu:
(\zmin\otimes_\B \zmin)\to \zmin$.  On simple tensors, $\mu$ is given by the equality
$
\mu(q^ a\otimes_\B q^b)=q^{a+b}
$
and thus,  since addition in $ \zmin$ corresponds to taking the $\inf$, one derives for arbitrary elements 
\begin{equation}\label{squaremap}
\mu\left(\sum q^{n_i}\otimes_\B q^{m_i}\right)=q^\alpha, \   \alpha= \inf (n_i+m_i).
\end{equation}
\begin{figure}
\begin{center}
\includegraphics[scale=0.4]{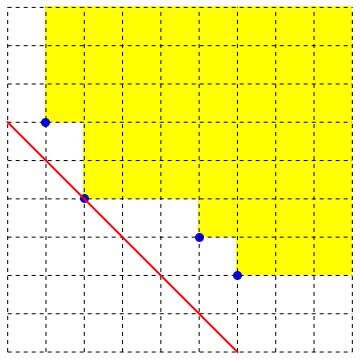}
\caption{Diagonal in the square, the morphism $\mu$.}\label{adsub1}
\end{center}
\end{figure}
\no Guided by this example of the product, we shall define the Frobenius correspondences $\cC_r$ for each positive rational number $r\in \Q_+^\times$ using the composition $\mu\circ \fr_{n,m}$ for $r=n/m$. 

\begin{prop}\label{frobcor} 
$(i)$~The range of the composite $\mu\circ \fr_{n,m}(\zmin^+\otimes_\B \zmin^+)\subset \zminp$ depends, up to canonical isomorphism, only upon the ratio $r=n/m$. Assuming $(n,m)=1$, this range contains the ideal
$
\{q^a\mid a\geq (n-1)(m-1)\}\subset \zminp
$.\newline
$(ii)$~The range of $\mu\circ \fr_{n,m}(\zmin^+\otimes_\B \zmin^+)\subset \zmin^+$ ($(n,m)=1$) is the semiring $F(n,m)$ generated by two elements $X,Y$ such that $X^m=Y^n$ and where the addition comes from a total order. If  $(n,m)=1$  and $n,m\neq 1$, the subset  $\{n,m\}\subset \N$ of such pairs is determined by the semiring $F(n,m)$.\newline
$(iii)$~Let $r=n/m$, $q\in(0,1)$  and let 
$$
m_r: \zmin\otimes_\B \zmin\to \rmax,\quad  m_r\left(\sum (q^{n_i}\otimes_\B q^{m_i})\right)=q^\alpha, \   \alpha= \inf (r n_i+ m_i).
$$
Up to canonical isomorphism of their ranges, the  morphisms $\mu\circ \fr_{n,m}$ and $m_r$ are equal.
\end{prop}
\proof 
$(i)$~The range $R=\mu\circ \fr_{n,m}(\zmin^+\otimes_\B \zmin^+)\subset \zminp$ consists of the $q^{na+mb}$ where $a,b\in \N$. Let $R_{n,m}=\{na+mb\mid a,b\in \N\}$. Let us first assume that $n$ and $m$ are relatively prime. Then as for the simplest case of the Frobenius problem, if one assumes that $(n,m)=1$, one has by a result of Sylvester
$
x\in R_{n,m} \  \forall x\geq (n-1)(m-1), \    (n-1)(m-1)-1\notin R_{n,m}.
$
When $n,m$ are arbitrary, let $d$ be their GCD, then $R_{n,m}=dR_{n',m'}$ with $n=dn'$ and $m=dm'$ and one  gets $(i)$.\newline
$(ii)$~One can compare the range $R=\mu\circ \fr_{n,m}(\zmin^+\otimes_\B \zmin^+)\subset \zminp$ with the semiring $F$ generated by two elements $X,Y$ such that $X^m=Y^n$ and such that the addition comes from a total order. One has a homomorphism $\rho:F\to R$ given by $\rho(X):=q^n$ and $\rho(Y):=q^m$. This morphism is surjective by construction. Any element of $F$ is uniquely of the form $X^aY^b$ where $b\in \{0,\ldots n-1\}$. One has 
$
\rho(X^aY^b)=q^{na}q^{mb}=q^{na+mb}
$, 
where $b\in \{0,\ldots n-1\}$ is uniquely determined from $z=na +mb$ by the equality $z\equiv mb$ mod. $n$ and the fact that $m$ is prime to $n$. This shows that $\rho$ is injective and the hypothesis that the addition in $F$ comes from a total order shows that $\rho$ is an isomorphism. Let us now show how one recovers the ratio $n/m$ or its inverse from the semiring $F$. In fact, we only use the multiplicative mono\"{\i}d of $F$. We assume that $n$ and $m$ are relatively prime and neither of them is equal to $1$. Then, an equality of the form $n=na+mb$ with $a$ and $b$ non-negative integers, implies that $a=1$ and $m=0$. Indeed, the only other possibility is $a=0$ but then one gets $n=mb$ which is not possible. This shows, using the isomorphism $\rho$ that the elements $X,Y$ of $F$ are the only ones which are indecomposable for the product. It follows that the two integers $n,m$ are uniquely determined from the relation $X^m=Y^n$, with $m$ prime to $n$. \newline 
$(iii)$~By \eqref{squaremap} one has
\begin{equation}\label{squaremap1}
\mu\circ \fr_{n,m}\left(\sum q^{n_i}\otimes_\B q^{m_i}\right)=q^\alpha, \quad   \alpha= \inf (nn_i+mm_i).
\end{equation}
Since $ \inf (nn_i+mm_i)=m \inf (r n_i+ m_i)$ with $r=n/m$, one gets the required equality with $m_r$. \endproof

The statement $(iii)$ of Proposition \ref{frobcor} shows how to define the Frobenius correspondence associated to a positive real number 
\begin{figure}
\begin{center}
\includegraphics[scale=0.4]{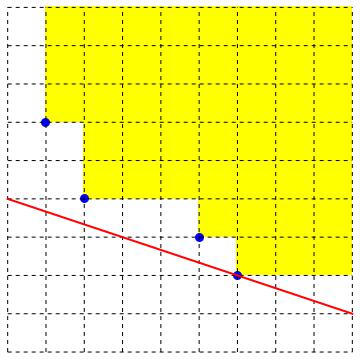}
\caption{The morphism $m_r$.}\label{adsub1}
\end{center}
\end{figure}
\begin{prop}\label{frobcorres} 
$(i)$~Let $\lambda\in \R_+^*$ and $q\in (0,1)$, then the following defines a homomorphism of semirings 
\begin{equation}\label{frobc1}
\cF(\lambda,q):\zmin\otimes_\B \zmin\to\rmax,\qquad \cF(\lambda,q)\left(\sum q^{n_i}\otimes_\B q^{m_i}\right)=q^\alpha, \   \alpha= \inf (\lambda n_i+m_i).
\end{equation} 
$(ii)$~The  semiring $\cR(\lambda):=\cF(\lambda,q)( \zmin^+\otimes_\B \zmin^+)$  does not depend, up to canonical isomorphism, of the choice of $q\in (0,1)$.\newline
$(iii)$~The semirings $\cR(\lambda)$ and $\cR(\lambda')$, for $\lambda,\lambda'\notin\Q_+^* $, are isomorphic if and only if $\lambda'=\lambda$ or $\lambda'=1/\lambda$.
\end{prop}
\proof $(i)$~We use the identification of $\zmin\otimes_\B \zmin$ with $\sub(\Z\times\Z)$ for which the two operations are given by the union $E\cup E'$ and the sum $E+ E'$.  The expression of $\cF(\lambda,q)$ in terms of $E\in \sub(\Z\times\Z)$ is, using the linear form $L_\lambda: \R^2\to \R$, $L_\lambda(a,b):=\lambda a+b$
$$
\cF(\lambda,q)(E)=q^\alpha, \  \  \alpha= \inf_{x\in E} L_\lambda(x).
$$
One has 
$$
 \inf_{E\cup E'} L_\lambda=\inf \{ \inf_{E} L_\lambda, \inf_{E'} L_\lambda\}, \  \   \inf_{E+ E'} L_\lambda=
\inf_{E} L_\lambda + \inf_{E'} L_\lambda
$$
which shows that $\cF(\lambda,q)$ is a homomorphism.\newline
$(ii)$~Let $q'\in (0,1)$, then there exists $u\in \R_+^\times$ such that $q'=q^u$. It then follows that  
$
\cF(\lambda,q')=\fr_u\circ \cF(\lambda,q)
$, 
where $\fr_u\in \Aut(\rmax)$, $\fr_u(x)=x^u$, $\forall x\in \rmax$.\newline
$(iii)$~The  semiring $\cR(\lambda)$ is a mono\"{\i}d for the multiplication and a totally ordered set using its additive structure. As a mono\"{\i}d, it contains two elements $X,Y$ which are indecomposable as products (and are $\neq 1$). We first choose $X=q^\lambda$ and $Y=q$. The total ordering allows one to compare $X^a$ with $Y^b$ for any positive integers $a, b\in \N$. Let
$$
D_-:=\{b/a\mid X^a<Y^b\}, \qquad  D_+:=\{b/a\mid X^a>Y^b\}.
$$
One has $X^a<Y^b$ if and only if $\lambda a>b$, \ie $b/a<\lambda$ so that $D_-=\{b/a\mid b/a<\lambda\}$ and similarly $D_+=\{b/a\mid \lambda<b/a\}$. Thus one recovers from the semiring $\cR(\lambda)$ the Dedekind cut which defines $\lambda$. The choice of the other of the two generators replaces $\lambda$ by $1/\lambda$. \endproof

\begin{rem}\label{getlamback}{\rm Notice that the parameter $\lambda\in \R_+^\times$ is uniquely specified by the pair of semiring homomorphisms $\ell:\zmin\to R$, $r:\zmin\to R$, where $\ell(q^n)=q^{n\lambda}$ and $r(q^n)=q^n$, while $R$ is the semiring generated by $\ell(\zmin)r(\zmin)$ which does not depend upon the choice of $q\in (0,1)$. Indeed, as above one obtains the Dedekind cut which defines $\lambda$ by comparing the powers of $X=q^\lambda$ and $Y=q$.
}\end{rem}
The congruence 
\begin{equation}\label{morphism}
x\sim y \iff \cF(\lambda,q)(x)=\cF(\lambda,q)(y)
\end{equation}
associated to the  morphism 
$\cF(\lambda,q): \zmin\otimes_\B \zmin\to \rmax
$
is independent of $q$. It is important to note that this congruence 
can be obtained on any finite subset $Z\subset \zmin\otimes_\B \zmin$ by approximating  $\lambda$ by rational values $\lambda'\in \Q_+^\times$. Indeed, recall that a best rational approximation to a real number $\lambda$ is a rational number $n/d$, $d > 0$, that is closer to $\lambda$ than any approximation with a smaller or equal denominator. One has
\begin{lem}\label{ratapprox} Let $\lambda\in \R_+^\times$ and $n/d$ be a best rational approximation of $\lambda$. Then replacing $\lambda$ by $n/d$ does not change the following sets of rational numbers
$$
D_-(\lambda,d):=\{b/a\mid b/a<\lambda, \  a< d\}, \  \  D_+(\lambda,d):=\{b/a\mid b/a>\lambda, \  a< d\}.
$$
\end{lem}
\proof For $\lambda_1\leq \lambda_2$ one has $D_-(\lambda_1,d)\subset D_-(\lambda_2,d)$,~ $D_+(\lambda_1,d)\supset D_+(\lambda_2,d)$ and 
$$
D_-(\lambda_2,d)\setminus D_-(\lambda_1,d)=\{b/a\mid \lambda_1 \leq b/a<\lambda_2, \  a< d\}.
$$
Thus, if $D_-(\lambda_2,d)\neq D_-(\lambda_1,d)$, there exists $b/a\in [\lambda_1,\lambda_2)$ with $a< d$. Since the interval between  $\lambda$ and $n/d$ does not contain any rational number $b/a$ with $a<d$ one gets the equality $D_-(\lambda,d)=D_-(n/d,d)$. Similarly one has $D_+(\lambda,d)=D_+(n/d,d)$.\endproof

\begin{prop}\label{approxprop}
Let $\lambda\in \R_+^\times$ be irrational and $x,y\in \zmin\otimes_\B \zmin$. Then one has $\cF(\lambda,q)(x)=
\cF(\lambda,q)(y)$ if and only if for all sufficiently good approximations~ $n/m\sim \lambda$ one has 
$$
\mu\circ \fr_{n,m}(x)=\mu\circ \fr_{n,m}(y)
$$
where $\mu:\zmin\otimes_\B \zmin\to \zmin$ is the semiring multiplication. 
\end{prop}
\proof $(i)$~We fix $q\in (0,1)$. After multiplying $x$ and $y$ by $q^k\otimes  q^k$ for $k$ large enough, we can assume that $x,y\in \nboplus$. Let $x=\sum q^{n_i}\otimes_\B q^{m_i}$, $y=\sum q^{r_i}\otimes_\B q^{s_i}$, with $n_i,m_i,r_i,s_i\in \N$. Then by \eqref{frobc1} one has $$
\cF(\lambda,q)(x)=
\cF(\lambda,q)(y) \iff \inf (\lambda n_i+m_i)=\inf (\lambda r_i+s_i)
$$
Let then $n/d$ be a best rational approximation of $\lambda$ with $d>\sup\{n_i,r_i  \}$. We show that the equality 
$\inf (\lambda n_i+m_i)=\inf (\lambda r_i+s_i)$ is unaffected if one replaces $\lambda$ by $n/d$. To check this, it is enough to show that for any positive integers $a,a',b,b'$ with $a<d$, $a'<d$, one has 
$$
a\lambda +b<a'\lambda +b'\iff an/d +b<a' n/d +b'.
$$
This amounts to show that the comparison of $(a-a')\lambda$ with $b'-b$ is unaffected when one replaces $\lambda$ by $n/d$ which follows from Lemma \ref{ratapprox}.\endproof

Proposition \ref{approxprop} shows in a precise manner that for each $\lambda\in \R_+^\times$, not in $\Q_+^\times$, the congruence relation $\cC_\lambda$  on $\nbo$ defined by 
\eqref{morphism} 
is in fact the limit of the $\cC_r$, as the rational $r\to \lambda$ and in turns the congruence $\cC_r$ for $r=n/m$ is simply given by composing the endomorphism $\fr_{n,m}$ with the diagonal relation $\cC_1$ which corresponds to the product $\mu$ in the semiring $\zmin$. Thus one can write 
\begin{equation}\label{congruence1}
  \cC_\lambda= \lim_{n/m\to \lambda}\cC_1\circ \fr_{n,m}.
\end{equation}

\begin{defn} We let $\conv(\Z\times \Z)$ be the set of closed, convex subsets 
$$
C\subset \R^2, \quad C+\R_+^2=C, \quad \exists z\in \R^2, ~C\subset z+\R_+^2
$$ such that the extreme points $\partial C$ belong to $\Z\times \Z$.
\end{defn}

We endow $\conv(\Z\times \Z)$ with the following two operations, the first is the convex hull of the union
\begin{equation}\label{consum}
{\rm conv}(C\cup C')=\{\alpha x+(1-\alpha)x'\mid \alpha \in [0,1],\  \  x\in C, \  x'\in C'\}
\end{equation}
and the second is the sum $C+C'$.

\begin{lem}\label{approx0} Endowed with the above two operations $\conv(\Z\times \Z)$ is a semiring of characteristic $1$.
\end{lem}
\proof First notice that the operation defined in \eqref{consum} is commutative, associative, and admits $\emptyset$ as the neutral element. Similarly, the second operation $C+C'$  is commutative, associative, and admits $Q=\R_+^2$ as neutral element. Moreover, for three elements $C,C',C''\in \conv(\Z\times \Z)$, one has 
$$
{\rm conv}(C\cup C')+C''=\{\alpha x+(1-\alpha)x'+x''\mid \alpha \in [0,1],\   x\in C, \  x'\in C', \ x''\in C''\}
$$
which coincides, using the convexity of $C''$, with 
$$
{\rm conv}((C+C'')\cup (C'+C''))=\{\alpha (x+y)+(1-\alpha)(x'+z)\mid \alpha \in [0,1],\   x\in C, \  x'\in C', \ y,z\in C''\}
$$
Finally, by construction, the first operation is evidently idempotent.\endproof

The role of the semiring $\conv(\Z\times \Z)$ will be clarified in Proposition \ref{conv} where we show that it gives the reduction of $\nbo$.

\begin{prop}\label{approx}
Let $x,y\in \zmin\otimes_\B \zmin$. One has:  $\cF(\lambda,q)(x)=
\cF(\lambda,q)(y)$ $\forall\lambda\in \R_+^\times$ if and only if $\mu\circ \fr_{n,m}(x)=\mu\circ \fr_{n,m}(y)$ for all pairs of positive integers $n,m\in \N$. This defines a congruence relation on the semiring $\zmin\otimes_\B \zmin$. The quotient semiring is the semiring $\conv(\Z\times \Z)$.
\end{prop}
\proof 
 The first part follows from Proposition \ref{approxprop}. For each $\lambda\in \R_+^\times$ the relation $\cF(\lambda,q)(x)=
\cF(\lambda,q)(y)$ is a congruence relation on $\nbo$ independent of the choice of $q\in (0,1)$. The conjunction of congruences is a congruence and thus it only remains to show that the quotient semiring is $\conv(\Z\times \Z)$. This follows since the convex closure of any $E\in \sub(\Z\times \Z)$ is obtained as the intersection of the half planes of the form
$$
H_{\lambda,\alpha}:=\{(x,y)\in \R^2\mid \lambda x+y\geq \alpha\}, \   \   \lambda\in \R_+^\times, \  \alpha \geq 0.
$$
which contain $E$.\endproof

\subsection{The reduced square $\arithc$}

Next, we investigate the relation between the semirings $\nbo$ and $\conv(\Z\times \Z)$. One has by construction a homomorphism of semirings $\gamma:\nbo\to\conv(\Z\times \Z)$ which in terms of the description of $\nbo$ as $\sub(\Z\times \Z)$ is simply given by $E\mapsto \gamma(E)$, where $\gamma(E)$ is the convex hull of $E$. Observe that while the operation $E\mapsto nE$ defines the endomorphism $ \fr_{n,n}$ of $\sub(\Z\times \Z)=\nbo$, it does not coincide with the operation $x\mapsto x^n$. Indeed, the latter map corresponds to $E\mapsto E+E+\ldots +E$ which differs in general from $nE$. In particular  the operation $x\mapsto x^n$ is not an endomorphism of $\nbo$ and it fails to be additive. The simplest example of this failure of additivity is the following ($n=2$)
$$
(q\otimes_\B 1+1\otimes_\B q)^2=q^2\otimes_\B 1+q\otimes_\B q +1\otimes_\B q^2\neq q^2\otimes_\B 1+1\otimes_\B q^2.
$$
By Proposition 4.43 of \cite{Golan}, the map $x\mapsto x^n$ is an injective endomorphism for any multiplicatively cancellative semiring of characteristic $1$. We thus conclude that $\nbo$ fails to be multiplicatively cancellative. In fact one checks that the equality $ac=bc$ for $a=(q\otimes_\B 1+1\otimes_\B q)^2$, $b=q^2\otimes_\B 1+1\otimes_\B q^2$, $c=q\otimes_\B 1+1\otimes_\B q$, while $a\neq b$ and $c\neq 0$.

We now show that the homomorphism $\gamma:\nbo\to\conv(\Z\times \Z)$ is the same as the homomorphism from $\nbo$ to its semiring of quotients, as defined in Chapter 11, Proposition 11.5 of \cite{Golan}. 

We use the notation Ceiling$(x):=\inf \{n\in\Z\mid n\geq x\}$.

\begin{lem}\label{lemconv} Let $a,b\in \N$, both $\neq 0$. Define an element $\sigma(a,b)\in \nboplus$ by
$$
\sigma(a,b)=\sum_{0\leq j\leq a} q^{a-j}\otimes_\B q^{{\rm Ceiling}(\frac{bj}{a})}.
$$
One has $\cF(b/a,q)(\sigma(a,b))=q^b$ and 
\begin{equation}\label{cancellation0}
x\in \nboplus \ \text{and}\  \cF(b/a,q)(x)=q^b \implies x+\sigma(a,b)=\sigma(a,b).
\end{equation}
Moreover, the following equality holds in $\nboplus$
\begin{equation}\label{cancellation1}
(q^a\otimes_\B 1+1\otimes_\B q^b)\sigma(a,b)=\sigma(a,b)\sigma(a,b).
\end{equation}
\end{lem}
\proof The equality $\cF(b/a,q)(\sigma(a,b))=q^b$ follows from
$$
\inf_{0\leq j\leq a}((b/a)(a-j)+{\rm Ceiling}(\frac{bj}{a})=b.
$$
Let $E\in \sub(\N\times \N)$ be the subset of $\N\times \N$ associated to $\sigma(a,b)\in \nboplus$. By construction one has 
$
E=\{(n,m)\in \N\times \N\mid (b/a) n+m\geq b\}.
$
Let $F\in \sub(\N\times \N)$ be the subset of $\N\times \N$ associated to $x\in \nboplus$. Then, if  $\cF(b/a,q)(x)=q^b$ one has $F\subset E$ and  since addition in $\sub(\N\times \N)$ is given by union one gets $x+\sigma(a,b)=\sigma(a,b)$.\newline
To prove \eqref{cancellation1} we show that both sides are equal to $\sigma(2a,2b)$. For the left hand side one has
$$
(q^a\otimes_\B 1+1\otimes_\B q^b)\sigma(a,b)=
\sum_{0\leq j\leq a} q^{2a-j}\otimes_\B q^{{\rm Ceiling}(\frac{bj}{a})}+\sum_{0\leq i\leq a} q^{a-i}\otimes_\B q^{b+{\rm Ceiling}(\frac{bi}{a})}=
$$
$$
=\sum_{0\leq j\leq a} q^{2a-j}\otimes_\B q^{{\rm Ceiling}(\frac{2bj}{2a})}+\sum_{a\leq j\leq 2a} q^{2a-j}\otimes_\B q^{{\rm Ceiling}(\frac{2bj}{2a})}=\sigma(2a,2b)
$$
where we used the idempotent addition so that the common term only counts once. Next one has
$
(q^a\otimes_\B 1+1\otimes_\B q^b)+\sigma(a,b)=\sigma(a,b)
$
and thus for the right hand side of \eqref{cancellation1} one gets
$$
\sigma(a,b)\sigma(a,b)=\left((q^a\otimes_\B 1+1\otimes_\B q^b)+\sigma(a,b)     \right)\sigma(a,b)=\sigma(2a,2b)+\sigma(a,b)\sigma(a,b).
$$
Since $\cF(b/a,q)(\sigma(a,b))=q^b$ one has $\cF(2b/2a,q)(\sigma(a,b)\sigma(a,b))=q^{2b}$ and by applying \eqref{cancellation0} one derives 
$
\sigma(a,b)\sigma(a,b)+\sigma(2a,2b)=\sigma(2a,2b).
$
We thus get as required
$
(q^a\otimes_\B 1+1\otimes_\B q^b)\sigma(a,b)=\sigma(2a,2b)=\sigma(a,b)\sigma(a,b)
$.
\endproof
\begin{figure}
\begin{center}
\includegraphics[scale=0.4]{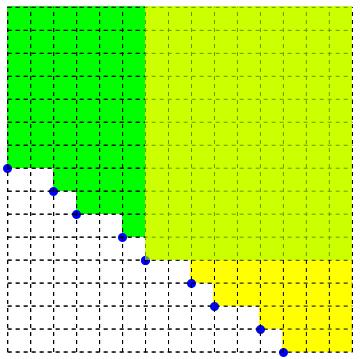}
\caption{The equality $(q^a\otimes_\B 1+1\otimes_\B q^b)\sigma(a,b)=\sigma(2a,2b)$
for $a=6$, $b=4$.}\label{adsub1}
\end{center}
\end{figure}
\begin{prop}\label{conv}
$(i)$~The semiring $\conv(\Z\times \Z)$ is multiplicatively cancellative.\newline
$(ii)$~The homomorphism $\gamma:\nbo\to\conv(\Z\times \Z)$ coincides with the homomorphism from $\nbo$ to its image in its semiring of quotients. \newline
$(iii)$~
Let $R$ be a multiplicatively cancellative semiring and $\rho:\nbo\to R$
a homomorphism such that $\rho^{-1}(\{0\})=\{0\}$. Then, there exists a unique semiring homomorphism $\rho':\conv(\Z\times \Z)\to R$ such that 
$\rho=\rho'\circ \gamma$.
\end{prop}
\proof $(i)$~An element $x\in\conv(\Z\times \Z)$ is uniquely specified by the elements $\cF(\lambda,q)(x)$ $\forall\lambda \in \R_+^\times$ (one can fix $q\in (0,1)$). Moreover if $x\neq 0$ one has $\cF(\lambda,q)(x)\neq 0$, $\forall \lambda \in \R_+^\times$. Since the semifield $\rmax$ is multiplicatively cancellative  one has for $a,b,c\in\conv(\N\times \N)$, $a\neq 0$
$$
ab=ac\implies \cF(\lambda,q)(ab)= \cF(\lambda,q)(ac)~ \forall\lambda \in \R_+^\times\implies 
 \cF(\lambda,q)(b)= \cF(\lambda,q)(c)~\forall\lambda \in \R_+^\times\implies b=c
$$
$(ii)$~Since we have shown that $\conv(\Z\times \Z)$ is multiplicatively cancellative, and moreover the homomorphism $\gamma:\nbo\to\conv(\Z\times \Z)$ fulfills $\gamma^{-1}(\{0\})=\{0\}$, we get that 
$$
 \exists c\in \nbo, \  c\neq 0, \  \  ca=cb \implies \gamma(a)=\gamma(b).
$$
To prove $(ii)$  it is enough to prove the converse, \ie to show that for any $a,a'\in\nbo$ one has
$$
\gamma(a)=\gamma(a')\implies \exists c\in \nbo, \  c\neq 0, \  \  ca=ca'.
$$
Multiplying $a$ and $a'$ by $q^n\otimes  q^n$, for $n$ large enough, we can assume that $a,a'\in \nboplus$. 
Let then $E,E'\in \sub(\N\times \N)$ be associated to $a,a'\in \nboplus$. The equality $\gamma(a)=\gamma(a')$ means that the convex hulls of $E$ and $E'$ are the same. In particular the extreme points are the same, \ie 
$
\partial \conv(E)=\partial \conv(E').
$
We denote this finite subset of $\N\times \N$ by $\cE=\{(a_i,b_i)\mid 1\leq i\leq n\}$, where the sequence $a_i$ is strictly increasing and $b_i$ is strictly decreasing. If $n=1$ we have necessarily $a=q^{a_1}\otimes_\B q^{b_1}$ and the same holds for $a'$ so that $a=a'$. We thus assume $n>1$ and let 
$$
x=\sum_{1\leq i\leq n}q^{a_i}\otimes_\B q^{b_i}, \  \  y=\sum_{1\leq i\leq n-1}(q^{a_i}\otimes_\B q^{b_{i+1}})\sigma(a_{i+1}-a_i,b_i-b_{i+1}).
$$
One has 
$$
x=\sum_{1\leq i\leq n-1}(q^{a_i}\otimes_\B q^{b_{i+1}})\left(q^{a_{i+1}-a_i}\otimes_\B 1+ 1 \otimes_\B q^{b_i-b_{i+1}}\right)
$$
and it follows from Lemma \ref{lemconv} that with $c:=\prod_{1\leq i\leq n-1} \sigma(a_{i+1}-a_i,b_i-b_{i+1})$ one has $xc=yc$ since 
$$
\left(q^{a_{i+1}-a_i}\otimes_\B 1+ 1 \otimes_\B q^{b_i-b_{i+1}}\right)c=\sigma(a_{i+1}-a_i,b_i-b_{i+1})c \qqq i, \ 1\leq i\leq n-1.
$$
Moreover for any $E\in \sub(\N\times \N)$ such that $\partial \conv(E)=\cE$ one has $\cE+Q\subset E\subset \conv(\cE)$. This implies that with $a\in \nboplus$ associated to $E\in \sub(\N\times \N)$, one has with the above notations, $a=a+x$ and $a+y=y$. Multiplying by $c$ and using $xc=yc$ one gets $ac=ac+xc=ac+yc=yc=xc$. Thus, since $\partial \conv(E)=\partial \conv(E')$, one gets $ac=a'c$ as required.\newline
$(iii)$~Since $R$ is multiplicatively cancellative, and   $\rho:\nbo\to R$ fulfills $\rho^{-1}(\{0\})=\{0\}$, we get that 
$$
 \exists c\in \nbo, \  c\neq 0, \  \  ca=cb \implies \rho(a)=\rho(b).
$$
Thus $\rho(a)$ only depends upon $\gamma(a)\in \conv(\Z\times \Z)$  and one obtains the required factorization. \endproof

\begin{defn}\label{site3} The {\em reduced square} $\arithc$ of the arithmetic site  is the topos $\wntb$ with the  {\em structure sheaf} $\conv(\Z\times \Z)$, viewed as a semiring in the topos.
\end{defn}

\begin{rem}\label{casereduced}{\rm One can prove in the same way  the analogue of Proposition  \ref{tensjustifsquare}, \ie that the set of points of the reduced square $\arithc$ over $\rmax$ is canonically isomorphic to the square of the set of points of $\arith$ over $\rmax$.
}\end{rem}

\section{Composition of  Frobenius correspondences}\label{sectcomp}

In this section we investigate the composition of  Frobenius correspondences. In order to relate directly the next results to those  announced in \cite{CCas}, we shall work out the tensor product using the restriction to $\zminp$ which simplifies the description of the reduced correspondences. The main result of this section is Theorem \ref{thmcomp} which describes the composition law of Frobenius correspondences. A subtle feature appears with respect to the expected composition rule $\Psi(\lambda)\circ \Psi(\lambda')=\Psi(\lambda\lambda')$. In fact, the equality holds when  $\lambda\lambda'\notin \Q$ but when one composes $\Psi(\lambda)\circ \Psi(\lambda')$  for irrational value of $\lambda$ and with $\lambda\lambda'\in \Q$, one obtains the tangential deformation of the correspondence $\Psi(\lambda\lambda')$. The reason behind this result is that for irrational $\lambda$, the Frobenius correspondence $\Psi(\lambda)$ has a flexibility which is automatically inherited by composition, and this property is not present for the Frobenius correspondence $\Psi(\lambda)$ when $\lambda$ is rational, but it is restored in the tangential deformation.

\subsection{Reduced correspondences}
 We view  the Frobenius correspondence associated to $\lambda \in \R_+^\times$ as the homomorphism $\cF(\lambda):\nboplus\to \cR(\lambda)$
defined in Proposition \ref{frobcorres}. By construction, the semiring $\cR(\lambda)$ is multiplicatively cancellative. Moreover the semiring $\nboplus$ is generated by the images $\iota_j(\zminp)$ of the two projections, \ie by the morphisms
$
\iota_1(q^n):=q^n\otimes_\B 1, \    \iota_2(q^n):=1\otimes_\B q^n \ \forall n\in \N.
$
Composing $\cF(\lambda)$ with the projections yields two morphisms $\zminp\to\cR(\lambda)$,  $\ell(\lambda):=\cF(\lambda)\circ \iota_1$ and $r(\lambda):=\cF(\lambda)\circ \iota_2$  which fulfill the following definition

\begin{defn}\label{defcorr} A {\em reduced correspondence} on the arithmetic site $\arith$ is given by a triple $(R,\ell,r)$, where $R$ is a \mc semiring, $\ell, r: \zminp\to R$ are semiring morphisms such that $\ell^{-1}(\{0\})=\{0\}$, $r^{-1}(\{0\})=\{0\}$
and  $R$ is generated by $\ell(\zminp)r(\zminp)$. 
\end{defn}

By construction, the Frobenius correspondence 
determines a reduced correspondence 
\begin{equation}\label{frobcordef}
\Psi(\lambda):=(R,\ell(\lambda),r(\lambda)), \  \  R:=\cR(\lambda), \  \  \ell(\lambda):=\cF(\lambda)\circ \iota_1, \  \  r(\lambda):=\cF(\lambda)\circ \iota_2.
\end{equation}
From \eqref{frobc1} one gets that the elements of $\cR(\lambda)$ are powers $q^\alpha$, where $\alpha \in \N+\lambda \N$ and  the morphisms $\ell(\lambda)$ and $r(\lambda)$ are described as follows
\begin{equation}\label{morphism3}
\ell(\lambda)(q^n)q^\alpha=q^{\alpha +n\lambda}, \    \  r(\lambda)(q^n)q^\alpha=q^{\alpha +n}.
\end{equation}

\begin{lem}\label{lemcor} Let $(R,\ell,r)$ be a reduced correspondence over $\arith$, then there exists a unique semiring homomorphism $\rho:\conv(\N\times \N)\to R$ such that $\ell=\rho\circ \gamma\circ \iota_1$ and $r=\rho\circ \gamma\circ \iota_2$.
\end{lem}
\proof By Proposition \ref{tensjustif} there exists a unique $\B$-linear map $\rho_0:\nboplus\to R$ such that 
$
\rho_0(q^a\otimes_\B q^b)=\ell(q^a)r(q^b)\ \forall a,b \in \N.
$
Moreover $\rho_0$ is multiplicative and hence it defines a homomorphism of semirings $\rho_0:\nboplus\to R$.  Let us show that $\rho_0^{-1}(\{0\})=\{0\}$. We recall that a semiring of characteristic $1$ is zero sum free
(\!\!\cite{Golan} page 150) \ie it satisfies the property 
\begin{equation}\label{zerosumfree}
x_j\in R,  \ \   \sum x_j=0\implies x_j=0,\  \forall j.
\end{equation}
Indeed, if $a+b=0$ then $a=a+a+b=a+b=0$. As $R$ is \mc by hypothesis, it has no zero divisors and $\ell(q^a)r(q^b)\neq 0\  \forall a,b \in \N$. This fact together with \eqref{zerosumfree} shows that $\rho_0^{-1}(\{0\})=\{0\}$.
We now apply Proposition \ref{conv} $(iii)$ and get that there exists a unique semiring homomorphism $\rho:\conv(\N\times \N)\to R$ such that $\rho_0=\rho\circ \gamma$. It follows  that $\ell=\rho\circ \gamma\circ \iota_1$ and $r=\rho\circ \gamma\circ \iota_2$. Since $\gamma:\nboplus\to\conv(\N\times \N)$ is surjective one gets the uniqueness of the homomorphism $\rho$. \endproof

\subsection{$\cR(\lambda)\otimes_{\zminp} \cR(\lambda')$, for $\lambda\lambda'\notin \Q\lambda'+\Q$}
We compute the composition of two Frobenius correspondences, viewed as reduced correspondences, \ie we consider the reduced correspondence associated to the formal expression
$$
\left(\cR(\lambda)\otimes_{\zminp} \cR(\lambda'), \   \ell(\lambda)\otimes \id, \  \id\otimes r(\lambda')\right).
$$
Reducing this expression means that one looks for the \mc semiring associated to the semiring  $\cR(\lambda)\otimes_{\zminp} \cR(\lambda')$ and then the sub-semiring generated by the images of  $\ell(\lambda)\otimes \id$ and $\id\otimes r(\lambda')$. In particular, we first study the meaning of the semiring $\cR(\lambda)\otimes_{\zminp} \cR(\lambda')$ and the associated \mc semiring. For $\lambda, \lambda'\in \R_+^\times$ and $q\in (0,1)$, we let $\cR(\lambda,\lambda')$ be the sub-semiring of $\rmax$ of elements of the form $q^\alpha$ for $\alpha \in \N+\lambda\N+\lambda'\N$. Up to canonical isomorphism, $\cR(\lambda,\lambda')$ is independent of the choice of $q\in (0,1)$. We denote as above  $\fr_u\in \Aut(\rmax)$, $\fr_u(x)=x^u$, $\forall x\in \rmax$. 
\begin{prop}\label{frobcorres1} 
$(i)$~Consider the map  
$$
\psi:\cR(\lambda)\times\cR(\lambda')\to \cR(\lambda\lambda',\lambda'),\quad \psi(a,b):=\fr_{\lambda'}(a)b\qqq a\in \cR(\lambda),\ b\in \cR(\lambda').
$$
Then $\psi$ is bilinear: $\psi(aa',bb')=\psi(a,b)\psi(a',b')$, $\forall a,a' \in \cR(\lambda)$, $b,b' \in \cR(\lambda')$ and it satisfies the equality
$$
\psi(r(\lambda)(x)a,b)=\psi(a,\ell(\lambda')(x)b)\qqq a\in \cR(\lambda),\ b\in \cR(\lambda'), \ x\in \zminp.
$$
$(ii)$~Let $R$ be a \mc semiring and $\phi:\cR(\lambda)\times\cR(\lambda')\to R$ be a bilinear map such that
$\phi(aa',bb')=\phi(a,b)\phi(a',b')$, $\forall a,a' \in \cR(\lambda)$, $b,b' \in \cR(\lambda')$ and
\begin{equation}\label{hypo1}
\phi(r(\lambda)(x)a,b)=\phi(a,\ell(\lambda')(x)b)\qqq a\in \cR(\lambda),\ b\in \cR(\lambda'), \ x\in \zminp.
\end{equation}
Then, assuming  $\lambda\lambda'\notin \Q\lambda'+\Q$,  there exists a unique homomorphism $\rho: \cR(\lambda\lambda',\lambda')\to R$ such that $\phi=\rho\circ\psi$.
\end{prop}
\proof $(i)$~By construction $\psi(a,b)$ is the product of the two homomorphisms $a\mapsto \fr_{\lambda'}(a)$
and $b\mapsto b$ hence the bilinearity and multiplicativity properties follow. Moreover, for $a=q^\alpha\in \cR(\lambda)$, $b=q^\beta\in \cR(\lambda')$ and $x= q^n\in \zminp$ one has, using \eqref{morphism3}
$$
\psi(r(\lambda)(x)a,b)=\psi(r(\lambda)(q^n)q^\alpha,q^\beta)=\psi(q^{\alpha+n},q^\beta)=q^{\lambda'(\alpha+n)+\beta}
$$
$$
\psi(a,\ell(\lambda')(x)b)=\psi(q^\alpha,\ell(\lambda')(q^n)q^\beta))=\psi(q^\alpha,q^{\beta+n\lambda'})=q^{\lambda'(\alpha+n)+\beta}.
$$
$(ii)$~Let $\alpha_i\in \N+\lambda \N, \   \beta_i\in \N+\lambda' \N$ be such that $\lambda'\alpha_1+\beta_1<\lambda'\alpha_2+\beta_2$. Let us show that 
\begin{equation}\label{simprule2}
\phi(q^{\alpha_1}, q^{\beta_1})+\phi(q^{\alpha_2}, q^{\beta_2})=\phi(q^{\alpha_1}, q^{\beta_1}).
\end{equation}
Since $R$ is multiplicatively cancellative, the map $R\ni x\mapsto x^n\in R$ is an injective endomorphism. Thus  it is enough to show that for some $n\in \N$ one has 
\begin{equation}\label{simprule3}
\phi(q^{\alpha_1}, q^{\beta_1})^n+\phi(q^{\alpha_2}, q^{\beta_2})^n=\phi(q^{\alpha_1}, q^{\beta_1})^n.
\end{equation}
Using the multiplicativity of $\phi$ one has $\phi(q^{\alpha_j}, q^{\beta_j})^n=\phi(q^{n\alpha_j}, q^{n\beta_j})$. Thus it is enough to show that for some $n\in \N$ one has 
\begin{equation}\label{simprule4}
\phi(q^{n\alpha_1}, q^{n\beta_1}) +\phi(q^{n\alpha_2}, q^{n\beta_2}) =\phi(q^{n\alpha_1}, q^{n\beta_1})
\end{equation}
By \eqref{hypo1} one has 
\begin{equation}\label{hypo2}
\phi(q^{\alpha}, q^{\beta+\lambda' k})=\phi(q^{\alpha+ k}, q^{\beta}).
\end{equation}
We can find $n\in\N$  and $k,k'\in \N$ such that 
$$
\alpha_1 \leq \frac kn, \  \frac{ k'}{n}\leq \alpha_2, \ \lambda'\frac kn+\beta_1< \lambda'\frac{ k'}{n}+\beta_2.
$$
One then has, using $n\alpha_1\leq k$ and the bilinearity of $\phi$
$$
\phi(q^{n\alpha_1}, q^{n\beta_1})=\phi(q^{n\alpha_1}, q^{n\beta_1})+\phi(q^{k}, q^{n\beta_1}).
$$
By \eqref{hypo2}, one has  $\phi(q^{k}, q^{n\beta_1})=\phi(1,q^{n\beta_1+\lambda'k})$ and since $\lambda'\frac kn+\beta_1< \lambda'\frac{ k'}{n}+\beta_2$ one gets
$$
\phi(1,q^{n\beta_1+\lambda'k})+\phi(1,q^{n\beta_2+\lambda'k'})=\phi(1,q^{n\beta_1+\lambda'k})
$$
so that 
$$
\phi(q^{k}, q^{n\beta_1})+\phi(q^{k'}, q^{n\beta_2})=\phi(q^{k}, q^{n\beta_1}).
$$
Moreover, using $k'\leq n\alpha_2$ one has 
$$
\phi(q^{k'}, q^{n\beta_2})=\phi(q^{k'}, q^{n\beta_2})+\phi(q^{n\alpha_2}, q^{n\beta_2}).
$$
Combining the above equalities gives
$$
\phi(q^{n\alpha_1}, q^{n\beta_1}) +\phi(q^{n\alpha_2}, q^{n\beta_2})=\phi(q^{n\alpha_1}, q^{n\beta_1})+\phi(q^{k}, q^{n\beta_1})+\phi(q^{n\alpha_2}, q^{n\beta_2})=
$$
$$
=\phi(q^{n\alpha_1}, q^{n\beta_1})+\phi(q^{k}, q^{n\beta_1})+\phi(q^{k'}, q^{n\beta_2})+\phi(q^{n\alpha_2}, q^{n\beta_2})
=\phi(q^{n\alpha_1}, q^{n\beta_1})+\phi(q^{k}, q^{n\beta_1})+\phi(q^{k'}, q^{n\beta_2})=
$$
$$
=\phi(q^{n\alpha_1}, q^{n\beta_1})+\phi(q^{k}, q^{n\beta_1})=\phi(q^{n\alpha_1}, q^{n\beta_1}).
$$
Thus we have shown \eqref{simprule4} and hence \eqref{simprule2}. Let us now assume that  $\lambda\lambda'\notin \Q\lambda'+\Q$. Every element of  $\cR(\lambda\lambda',\lambda')$ is uniquely of the form $q^z$, $z=x\lambda\lambda'+\beta$, $x\in \N$, $\beta\in \N\lambda'+\N$.  Define the map $\rho: \cR(\lambda\lambda',\lambda')\to R$ by $\rho(q^z):=\phi(q^{x\lambda},q^{\beta})$. It is multiplicative by construction. We show that it is additive, \ie that for $z,z'\in \N\lambda\lambda'+\N\lambda'+\N$, $z<z'$ one has $\rho(q^z)=\rho(q^z)+\rho(q^{z'})$.
This follows from \eqref{simprule2} applied to $\alpha_1=x\lambda$, $\beta_1=\beta$, $\alpha_2=x'\lambda$, $\beta_2=\beta'$. Finally one has $\phi(a,b)=\rho(\psi(a,b))$ for $a=q^\alpha\in \cR(\lambda)$, $b=q^\beta\in \cR(\lambda')$ since with $\alpha=x\lambda +y$,  one has $\psi(a,b)=q^{x\lambda\lambda'+y\lambda'+\beta}$ and
$$
\rho(\psi(a,b))=\phi(q^{x\lambda},q^{y\lambda'+\beta})=\phi(q^{x\lambda+y},q^{\beta})=\phi(a,b).
$$
\endproof

\subsection{$\cR(\lambda)\otimes_{\zminp} \cR(\lambda')$ for $\lambda\lambda'\in \Q\lambda'+\Q$}
In this subsection we consider the more delicate case when $\lambda$ and $\lambda'$ are irrational but $\lambda\lambda'\in \Q\lambda'+\Q$. We let $\germ$ be the semiring of  germs of continuous functions from a neighborhood of $0\in \R$ to $\rmax$, endowed with the pointwise operations.  We consider the sub-semiring $\cR_\epsilon(\lambda\lambda',\lambda')$ of  $\germ$ generated, for fixed $q\in (0,1)$, by $q$, $q^{\lambda'}$ and $\fr_{1+\epsilon}(q^{\lambda\lambda'})=q^{(1+\epsilon)\lambda\lambda'}$. Proposition \ref{frobcorres1} adapts as follows

\begin{prop}\label{frobcorres2} Let  $\lambda, \lambda'\in \R_+^\times$ be irrational and such that $\lambda\lambda'\in \Q\lambda'+\Q$.

$(i)$~Let $\psi:\cR(\lambda)\times\cR(\lambda')\to \cR_\epsilon(\lambda\lambda',\lambda')$ be given by 
$$
\psi(q^{\lambda i+j},b):=q^{(1+\epsilon)\lambda\lambda' i}q^{\lambda'j}\,b    \qqq q^{\lambda i+j}\in \cR(\lambda),\ b\in \cR(\lambda').
$$
Then $\psi$ is bilinear, $\psi(aa',bb')=\psi(a,b)\psi(a',b')$, $\forall a,a' \in \cR(\lambda)$, $b,b' \in \cR(\lambda')$ and 
\begin{equation}\label{compa}
\psi(r(\lambda)(x)a,b)=\psi(a,\ell(\lambda')(x)b)\qqq a\in \cR(\lambda),\ b\in \cR(\lambda'), \ x\in \zminp.
\end{equation}
$(ii)$~Let $R$ be a \mc semiring and $\phi:\cR(\lambda)\times\cR(\lambda')\to R$ be a bilinear map such that
$\phi(aa',bb')=\phi(a,b)\phi(a',b')$, $\forall a,a' \in \cR(\lambda)$, $b,b' \in \cR(\lambda')$ and
\begin{equation}\label{hypo1bis}
\phi(r(\lambda)(x)a,b)=\phi(a,\ell(\lambda')(x)b)\qqq a\in \cR(\lambda),\ b\in \cR(\lambda'), \ x\in \zminp.
\end{equation}
Then   there exists a unique homomorphism $\rho: \cR_\epsilon(\lambda\lambda',\lambda')\to R$ such that $\phi=\rho\circ\psi$.
\end{prop}
\proof $(i)$~Since $\lambda$ is irrational the map $\delta:\cR(\lambda)\to \germ$, $\delta(q^{\lambda i+j}):=q^{((1+\epsilon)\lambda i)+j}$ is an isomorphism with its image. Indeed, for $(i,j)\neq (i',j')$ one has a strict inequality of the form $\lambda i+j<\lambda i'+j'$ and this still holds for $\epsilon$ small enough when one replaces $\lambda$ by $(1+\epsilon)\lambda$.  One has 
$$
\psi(a,b)=\fr_{\lambda'}(\delta(a))\, b\qqq a\in \cR(\lambda),\ b\in \cR(\lambda')
$$
which shows that $\psi$ is bilinear and multiplicative. One checks directly the equality \eqref{compa}.\newline
$(ii)$~Let $U:=\phi(q^\lambda,1)$, $V:=\phi(q,1)=\phi(1,q^{\lambda'})$,  $W:=\phi(1,q)$. Let $U_\epsilon=\psi(q^\lambda,1)=q^{(1+\epsilon)\lambda\lambda'}$, $V_\epsilon=\psi(q,1)=\psi(1,q^{\lambda'})=q^{\lambda'}$,  $W_\epsilon=\psi(1,q)=q$. The equality $\phi=\rho\circ\psi$ implies that one must have 
\begin{equation}\label{unique}
\rho(U_\epsilon)=U, \    \rho(V_\epsilon)=V, \    \rho(W_\epsilon)=W  
\end{equation}
which shows the uniqueness of $\rho$, if it exists, since by construction $U_\epsilon,V_\epsilon,W_\epsilon$ generate $\cR_\epsilon(\lambda\lambda',\lambda')$. Let us  show the existence of $\rho$. As shown in the proof of Proposition \ref{frobcorres1}, one has the addition rule \eqref{simprule2} which means with our notations
$$
a,b,c,a',b',c'\in \N\ \text{and}\ \lambda\lambda' a+\lambda'b+c<\lambda\lambda' a'+\lambda'b'+c' \implies U^aV^bW^c+U^{a'}V^{b'}W^{c'}=U^aV^bW^c.
$$
The same rule holds by direct computation for $U_\epsilon,V_\epsilon,W_\epsilon$ instead of $U,V,W$. This shows that if one defines $\rho$ on monomials using multiplicativity and \eqref{unique}, the additivity follows provided one proves that for any $u>0$ 
$$
F(u)=\{(a,b,c)\in \N^3\mid \lambda\lambda' a+\lambda'b+c=u\}
$$
and for subsets $S,S'\subset F(u)$ one has
$$
\sum_S U_\epsilon^aV_\epsilon^bW_\epsilon^c=\sum_{S'} U_\epsilon^aV_\epsilon^bW_\epsilon^c \implies
\sum_S U^aV^bW^c=\sum_{S'} U^aV^bW^c.
$$
Let $s,t\in \Q$ be the unique rational numbers  such that $\lambda\lambda'=s\lambda'+t$. Let $m$ be the lcm of the denominators of $s$ and $t$. One has $s=f/m$, $t=g/m$ with $f,g\in\Z$. Let $(a,b,c)\in F(u)$, $(a',b',c')\in F(u)$, then  one gets $c'-c=(a-a')t$, $b'-b=(a-a')s$ and there exists $k\in \Z$ such that $a-a'=km$, $b'-b=kf$, $c'-c=kg$. Thus if $F(u)\neq \emptyset$ there exists $(a_0,b_0,c_0)\in F(u)$ and an integer $n\in \N$ such that 
$$
F(u)=\{(a_0,b_0,c_0)+j(m,-f,-g)\mid j\in \{0,\ldots,n\}\}.
$$
 Both $\cR_\epsilon(\lambda\lambda',\lambda')$ and $R$ are \mc semirings and thus they embed in the semifields of fractions ${\rm Frac}(\cR_\epsilon(\lambda\lambda',\lambda'))$ and ${\rm Frac}(R)$ {\em resp.} Let then
$$
Z_\epsilon:=U_\epsilon^mV_\epsilon^{-f}W_\epsilon^{-g}\in {\rm Frac}(\cR_\epsilon(\lambda\lambda',\lambda')),
\  \  Z:=U^mV^{-f}W^{-g}\in {\rm Frac}(R).
$$
What  remains to be shown is that for any subsets $S,S'\subset \{0,\ldots,n\}$ one has 
$$
\sum_S Z_\epsilon^j=\sum_{S'} Z_\epsilon^j \implies
\sum_S Z^j=\sum_{S'} Z^j.
$$
Now, $Z_\epsilon$ is given by the germ of the function $\epsilon\mapsto q^{(1+\epsilon)m\lambda\lambda'}q^{-f\lambda'-g}=q^{m\epsilon \lambda\lambda'}$ and one gets 
$$
\sum_S Z_\epsilon^j=\sum_{S'} Z_\epsilon^j\iff \inf(S)=\inf(S')\ \text{and} \ \max(S)=\max(S')
$$
Thus the proof of $(ii)$ follows from the following Lemma \ref{span}. \endproof

\begin{lem}\label{span} Let $F$ be a semifield of characteristic $1$ and $Z\in F$, $n\in \N$. Then the additive span of the $Z^i$ for $i\in \{0,\ldots, n\}$ is the set of $Z(i,j):=Z^i+Z^j$ for $0\leq i\leq j\leq n$.\newline
 For any subset $S\subset \{0,\ldots,n\}$ one has
\begin{equation}\label{compsum}
\sum_S Z_\epsilon^k=Z(i,j), \  i=\inf(S), \  j=\max(S).
\end{equation}
\end{lem}
\proof 
Let $s(i,j):=\sum_{i\leq x\leq j} Z^x$. One has the following analogue of Lemma \ref{lemconv}
\begin{equation}\label{cancellation3}
Z(i,j)s(i,j)=s(i,j)s(i,j).
\end{equation}
In fact, we first show
\begin{equation}\label{cancellation4}
Z(i,j)s(i,j)=s(2i,2j).
\end{equation}
One has
$$
Z(i,j)s(i,j)=\sum_{i\leq x\leq j} Z^{x+i}+\sum_{i\leq x\leq j} Z^{x+j}.
$$
In the first sum $x+i$ varies from $2i$ to $i+j$, while in the second $x+j$ varies from $i+j$ to $2j$ while the repetition at $i+j$ does not affect the sum since $1+1=1$. Thus we get \eqref{cancellation4}. Moreover one has 
\begin{equation}\label{cancellation5}
s(i,j)s(k,\ell)=s(i+k,j+\ell).
\end{equation}
Indeed, this amounts to check that 
$$
\{(x,n-x)+(y,m-y)\mid i\leq x\leq j, \  k\leq y\leq \ell\}=\{(z,n+m-z)\mid i+k\leq z\leq j+\ell\}.
$$
Now by \eqref{cancellation3} one has $Z(i,j)=s(i,j)$ and it follows that for any subset $S\subset \{0,\ldots, n\}$ one gets $\sum_S Z^k=Z(i,j)$ where $i=\inf S$ and $j=\max S$. This shows that the additive span of the $Z^i$ for $i\in \{0,\ldots, n\}$ is the set of $Z(i,j):=Z^i+Z^j$ for $0\leq i\leq j\leq n$ and that \eqref{compsum} holds.\endproof

\subsection{The composition $\Psi(\lambda)\circ \Psi(\lambda')$}

Let $\beps$ be the the sub-semiring of  $\germ$ generated, for fixed $q\in (0,1)$, by $q$ and $\theta_{1+\epsilon}(q)=q^{1+\epsilon}$. It is independent, up to canonical isomorphism, of the choice of $q\in (0,1)$.

\begin{defn}\label{defndef} The tangential deformation of the identity correspondence is given by the triple $(\beps,\ell_\epsilon,r_\epsilon)$ where  $\ell_\epsilon(q^n):=\theta_{1+\epsilon}(q^n)$ and $r(q^n):=q^n$, $\forall n\in \N$.
\end{defn}
We then obtain:
\begin{thm}\label{thmcomp} Let $\lambda, \lambda' \in \R_+^*$ such that $\lambda\lambda'\notin \Q$. The composition of the Frobenius correspondences is then given by 
$$
\Psi(\lambda)\circ \Psi(\lambda')=\Psi(\lambda\lambda')
$$
The same equality holds if $\lambda$ and $\lambda'$ are rational. When  $\lambda, \lambda'$ are irrational and $\lambda\lambda'\in \Q$,
$$
\Psi(\lambda)\circ \Psi(\lambda')=\id_\epsilon\circ \Psi(\lambda\lambda')
$$
where $\id_\epsilon$ is the tangential deformation of the identity correspondence.
\end{thm}
\proof Let us first assume that $\lambda\lambda'\notin \Q\lambda'+\Q$. By Proposition \ref{frobcorres1} the \mc reduction of $\cR(\lambda)\otimes_{\zminp} \cR(\lambda')$ is $\cR(\lambda\lambda',\lambda')$ and the left and right actions of $\zminp$ are given by 
$$
\ell(q^n)q^\alpha=q^{\lambda\lambda'n}q^\alpha, \  \  r(q^n)q^\alpha=q^{n+\alpha}
$$
Thus the sub-semiring of the \mc reduction of $\cR(\lambda)\otimes_{\zminp} \cR(\lambda')$ which is generated by $\ell(\zminp)$ and $r(\zminp)$ is $\cR(\lambda\lambda')$ and the left and right actions of $\zminp$ are the same as for the correspondence $\Psi(\lambda\lambda')$.
 Let then  $\lambda, \lambda'\in \R_+^*$ be irrational and such that $\lambda\lambda'\in \Q\lambda'+\Q$.  By Proposition \ref{frobcorres2} the \mc reduction of $\cR(\lambda)\otimes_{\zminp} \cR(\lambda')$ is $\cR_\epsilon(\lambda\lambda',\lambda')$ and the left and right actions of $\zminp$ are given by 
$$
\ell(q^n)X=q^{(1+\epsilon)\lambda\lambda'n}X, \  \  r(q^n)X=q^{n}X \qqq X\in \cR_\epsilon(\lambda\lambda',\lambda').
$$
Assume first that $\lambda\lambda'\notin \Q$. Let us show that the sub-semiring $R$ of  $\cR_\epsilon(\lambda\lambda',\lambda')$ generated by the $\ell(q^n)$ and $r(q^m)$ is isomorphic to $\cR(\lambda\lambda')$ and the left and right actions of $\zminp$ are the same as for the correspondence $\Psi(\lambda\lambda')$. Indeed $R$ is the sub-semiring of $\germ$ generated by the $q^{(1+\epsilon)\alpha n}$ and $q^m$ where $\alpha=\lambda\lambda'$ is irrational. Its elements are finite sums  of the  form 
$$
\sigma(\epsilon)=\sum q^{(1+\epsilon)\alpha n_j+m_j}
$$ 
where the sum is taken in $\germ$. Since $\alpha\notin \Q$ the real numbers $\alpha n_j+m_j$ are distinct and thus there exists a unique $j_0$ such that 
$$
\alpha n_{j_0}+m_{j_0}<\alpha n_j+m_j \qqq j\neq j_0
$$
This implies that for $\epsilon$ small enough one has 
$$
(1+\epsilon)\alpha n_{j_0}+m_{j_0}<(1+\epsilon)\alpha n_j+m_j \qqq j\neq j_0
$$
and thus one gets the  equality $\sigma(\epsilon)=q^{(1+\epsilon)\alpha n_{j_0}+m_{j_0}}$ in $\germ$. This shows that the evaluation at $\epsilon=0$ is an isomorphism $R\to \cR(\alpha)$.

 For $\lambda\lambda'=1$ the sub-semiring $R$ of  $\cR_\epsilon(\lambda\lambda',\lambda')$ generated by the $\ell(q^n)$ and $r(q^m)$ is isomorphic to $\beps$   with left and right actions given as in 
Definition \ref{defndef} and thus one gets 
$$
\Psi(\lambda)\circ \Psi(\lambda^{-1})= \id_\epsilon.
$$
Assume now that $\alpha=\lambda\lambda'\in \Q$. Then the sub-semiring $R$ of  $\cR_\epsilon(\lambda\lambda',\lambda')$ generated by the $\ell(q^n)$ and $r(q^m)$ is formed of finite sums  of the  form 
$$
\sigma(\epsilon)=\sum q^{(1+\epsilon)\alpha n_j+m_j}
$$ 
where the sum is taken in $\germ$. The left and right actions are given by 
\begin{equation}\label{expected}
\ell(q^n)X=q^{(1+\epsilon)\alpha n}X, \  \  r(q^n)X=q^{n}X \qqq X\in R.
\end{equation}
Let us compare this with the compositions $\id_\epsilon\circ \Psi(\alpha)$ and $\Psi(\alpha)\circ \id_\epsilon$. We first need to determine the \mc reduction of  $\beps\otimes_{\zminp} \cR(\alpha)$.

\begin{lem}\label{frobcorres4} Let  $\alpha\in \Q_+^*$.

$(i)$~Let $\psi:\beps\times\cR(\alpha)\to \germ$ be given by 
$$
\psi(q^{(1+\epsilon) i+j},b):=q^{(1+\epsilon)\alpha i}q^{\alpha j}b\qqq i,j\in \N,\ b\in \cR(\alpha)
$$
Then $\psi$ is bilinear, $\psi(aa',bb')=\psi(a,b)\psi(a',b')$, $\forall a,a' \in \beps$, $b,b' \in \cR(\alpha)$ and 
\begin{equation}\label{compa4bis}
\psi(r_\epsilon(x)a,b)=\psi(a,\ell(\alpha)(x)b)\qqq a\in \beps,\ b\in \cR(\alpha), \ x\in \zminp
\end{equation}
$(ii)$~Let $R$ be a \mc semiring and $\phi:\beps\times\cR(\alpha)\to R$ be a bilinear map such that
$\phi(aa',bb')=\phi(a,b)\phi(a',b')$, $\forall a,a' \in \beps$, $b,b' \in \cR(\alpha)$ and 
\begin{equation}\label{hypo4bis}
\phi(r_\epsilon(x)a,b)=\phi(a,\ell(\alpha)(x)b)\qqq a\in \beps,\ b\in \cR(\alpha), \ x\in \zminp
\end{equation}
Then   there exists a unique homomorphism $\rho: \cR_\epsilon(\alpha,\alpha)\to R$ such that $\phi=\rho\circ\psi$.
\end{lem}
\proof $(i)$~One has $\psi(a,b)=\fr_\alpha(a)b$ which gives the required properties. 

 $(ii)$~Let $U:=\phi(q^{(1+\epsilon)},1)$, $V:=\phi(q,1)=\phi(1,q^{\alpha})$,  $W:=\phi(1,q)$. Let us show that  :
$$
a,b,c,a',b',c'\in \N\ \&\ \alpha a+\alpha b+c<\alpha a'+\alpha b'+c' \implies U^aV^bW^c+U^{a'}V^{b'}W^{c'}=U^aV^bW^c
$$
Let $a_i\in \N, \   \beta_i\in \N+\alpha \N$ be such that $\alpha a_1+\beta_1<\alpha a_2+\beta_2$. Let us show that 
\begin{equation}\label{simprule2bis}
\phi(q^{(1+\epsilon)a_1}, q^{\beta_1})+\phi(q^{(1+\epsilon)a_2}, q^{\beta_2})=\phi(q^{(1+\epsilon)a_1}, q^{\beta_1})
\end{equation}
Since $R$ is \mc the map $R\ni x\mapsto x^n\in R$ is an injective endomorphism. Thus it is enough to show that for some $n\in \N$ one has 
\begin{equation}\label{simprule4bis}
\phi(q^{(1+\epsilon)n a_1}, q^{n\beta_1})+\phi(q^{(1+\epsilon)n a_2}, q^{n\beta_2})=\phi(q^{(1+\epsilon)n a_1}, q^{n\beta_1})
\end{equation}
We can find $n\in\N$  and $k,k'\in \N$ such that 
$$
a_1 < \frac kn, \  \frac{ k'}{n}< a_2, \ \alpha\frac kn+\beta_1< \alpha\frac{ k'}{n}+\beta_2
$$
One then has, using $n a_1< k$, that in $\beps$, $q^{(1+\epsilon)n a_1}+ q^k=q^{(1+\epsilon)n a_1}$ and 
$$
\phi(q^{(1+\epsilon)n a_1}, q^{n\beta_1})=\phi(q^{(1+\epsilon)n a_1}, q^{n\beta_1})+\phi(q^{k}, q^{n\beta_1})
$$
Using \eqref{hypo4bis}
 one has  $\phi(q^{k}, q^{n\beta_1})=\phi(1,q^{n\beta_1+\alpha k})$ and since $\alpha\frac kn+\beta_1< \alpha\frac{ k'}{n}+\beta_2$ one gets
$$
\phi(1,q^{n\beta_1+\alpha k})+\phi(1,q^{n\beta_2+\alpha k'})=\phi(1,q^{n\beta_1+\alpha k})
$$
so that 
$$
\phi(q^{k}, q^{n\beta_1})+\phi(q^{k'}, q^{n\beta_2})=\phi(q^{k}, q^{n\beta_1})
$$
Moreover, using $k'< n a_2$ one has $q^{(1+\epsilon)n a_2}+ q^{k'}=q^{k'}$ in $\beps$ so that
$$
\phi(q^{k'}, q^{n\beta_2})=\phi(q^{k'}, q^{n\beta_2})+\phi(q^{(1+\epsilon)n a_2}, q^{n\beta_2})
$$
Combining the above equalities gives
$$
\phi(q^{(1+\epsilon)n a_1}, q^{n\beta_1})+\phi(q^{(1+\epsilon)n a_2}, q^{n\beta_2})=\phi(q^{(1+\epsilon)n a_1}, q^{n\beta_1})+\phi(q^{k}, q^{n\beta_1})+\phi(q^{(1+\epsilon)n a_2}, q^{n\beta_2})
$$
$$
=\phi(q^{(1+\epsilon)n a_1}, q^{n\beta_1})+\phi(q^{k}, q^{n\beta_1})+\phi(q^{k'}, q^{n\beta_2})+\phi(q^{(1+\epsilon)n a_2}, q^{n\beta_2})$$ $$
=\phi(q^{(1+\epsilon)n a_1}, q^{n\beta_1})+\phi(q^{k}, q^{n\beta_1})+\phi(q^{k'}, q^{n\beta_2})
=\phi(q^{(1+\epsilon)n a_1}, q^{n\beta_1})+\phi(q^{k}, q^{n\beta_1})=\phi(q^{(1+\epsilon)n a_1}, q^{n\beta_1}).
$$

The sub-semiring of $\germ$ generated by the range of $\psi$ is generated by the $q^{(1+\epsilon)\alpha i+\alpha j+k}$ for $i,j,k\in\N$.
Let $U_\epsilon=\psi(q^{(1+\epsilon)},1)=q^{(1+\epsilon)\alpha}$, $V_\epsilon=\psi(q,1)=\psi(1,q^{\alpha})=q^{\alpha}$,  $W_\epsilon=\psi(1,q)=q$. The equality $\phi=\rho\circ\psi$ implies that one must have 
\begin{equation}\label{uniquebis}
\rho(U_\epsilon)=U, \    \rho(V_\epsilon)=V, \    \rho(W_\epsilon)=W, \    
\end{equation}
which shows the uniqueness of $\rho$, if it exists, since by construction $U_\epsilon,V_\epsilon,W_\epsilon$ generate $\cR_\epsilon(\alpha,\alpha)$. Let us show the existence of $\rho$. As shown above
$$
a,b,c,a',b',c'\in \N\ \&\ \alpha a+\alpha b+c<\alpha a'+\alpha b'+c' \implies U^aV^bW^c+U^{a'}V^{b'}W^{c'}=U^aV^bW^c
$$
and the same rule holds by direct computation for $U_\epsilon,V_\epsilon,W_\epsilon$ instead of $U,V,W$. This shows that if one defines $\rho$ on monomials using multiplicativity and \eqref{uniquebis}, the additivity follows provided one proves that for any $u>0$ 
$$
F(u)=\{(a,\beta)\in \N\times (\alpha\N+\N)\mid \alpha a+\beta=u\}
$$
and subsets $S,S'\subset F(u)$ one has
$$
\sum_S U_\epsilon^aW_\epsilon^\beta=\sum_{S'} U_\epsilon^aW_\epsilon^\beta \implies
\sum_S U^a W^\beta=\sum_{S'} U^a W^\beta
$$
Here we use the notations
$$
W_\epsilon^\beta:=V_\epsilon^{b} W_\epsilon^c, \  W^\beta:=V^{b} W^c, \qqq  \beta=\alpha b+c, \ b,c\in\N
$$
The products $V_\epsilon^{b} W_\epsilon^c$ and $V^{b} W^c$ only depend on $\beta=\alpha b+c$.
Now assume $F(u)\neq \emptyset$, then $u\in \alpha\N+\N$ by construction since $\alpha a+\beta=u$ for some $(a,\beta)\in \N\times (\alpha\N+\N)$. The pairs $(a,\beta)\in \N\times (\alpha\N+\N)$ such that $\alpha a+\beta=u$ thus contain $(0,u)$ and all the $(a,u-a\alpha)$ such that $a\in \N$ and $u-a\alpha\in \alpha\N+\N$. There is a largest $t\in \N$ such that $u-t\alpha\in \alpha\N+\N$ and one then has 
$$
F(u)=\{(a,u-a\alpha)\mid 0\leq  a \leq t\}
$$
 Both $\cR_\epsilon(\alpha,\alpha)$ and $R$ are \mc semirings and thus they embed in the semifields of fractions ${\rm Frac}(\cR_\epsilon(\alpha,\alpha))$ and ${\rm Frac}(R)$. Let then
$$
Z_\epsilon:=U_\epsilon V_\epsilon^{-1}\in {\rm Frac}(\cR_\epsilon(\alpha,\alpha)),
\  \  Z:=UV^{-1}\in {\rm Frac}(R).
$$
What  remains to be shown is that for any subsets $S,S'\subset \{0,\ldots,t\}$ one has 
$$
\sum_S Z_\epsilon^j=\sum_{S'} Z_\epsilon^j \implies
\sum_S Z^j=\sum_{S'} Z^j
$$
Now $Z_\epsilon$ is given by the germ of the function $\epsilon\mapsto q^{(1+\epsilon)\alpha}q^{-\alpha}=q^{\epsilon \alpha}$ and one gets 
$$
\sum_S Z_\epsilon^j=\sum_{S'} Z_\epsilon^j\iff \inf(S)=\inf(S')\ \& \ \max(S)=\max(S')
$$
Thus the proof of $(ii)$ follows from  Lemma \ref{span}. \endproof

In order to obtain the composition $\id_\epsilon\circ \Psi(\lambda\lambda')$ one then takes the left and right actions of $\zminp$ in the sub-semiring  of $\germ$ generated by the range of $\psi$ \ie generated by the $q^{(1+\epsilon)\alpha i+\alpha j+k}$ for $i,j,k\in\N$. One gets the sub-semiring  of $\germ$  generated by the $q^{(1+\epsilon)\alpha i +k}$ for $i,k\in\N$, with left action of $q$ given by $i\mapsto i+1$ and right action by $k\mapsto k+1$. This coincides with \eqref{expected} and thus shows that 
$$
\Psi(\lambda)\circ \Psi(\lambda')=\id_\epsilon\circ \Psi(\lambda\lambda')
$$
where $\id_\epsilon$ is the tangential deformation of the identity correspondence.

\section{The structure of the point in noncommutative geometry}\label{sectncg}

This section is motivated by the classification of matro\"ids discovered by J. Dixmier in \cite{Dixmier}. It turns out that this classification involves the same noncommutative space as the space of points of the topos $\wnt$. The conceptual reason behind this coincidence comes from the structure of the topological space reduced to a single point when considered from the point of view of noncommutative geometry. As we show below this structure leads naturally to the topos $\wnt$. \newline
In noncommutative geometry a topological space is encoded by a $C^*$-algebra. Moreover  the large class of examples coming from spaces of leaves of foliations made it clear that  Morita equivalent $C^*$-algebras represent the same underlying space. In each Morita equivalence class of (separable) $C^*$-algebras one finds a unique (up to isomorphism) representative $A$ which is {\em stable} \ie  such that $A$ is isomorphic with the tensor product $A\otimes \cK$ of $A$ by the $C^*$-algebra $\cK$ of compact operators in Hilbert space. In particular, the single point is represented by the $C^*$-algebra $\cK$. This $C^*$-algebra is the natural home for the theory of infinitesimals (\cf \cite{NCG}). All the (star) automorphisms of the $C^*$-algebra $\cK$ are inner, \ie they are implemented by unitaries  in the multiplier algebra. But this algebra admits non-trivial endomorphisms and in fact one has  the following well-known result  where an endomorphism is called non-degenerate if it transforms approximate units into approximate units.
\begin{thm}\label{pointstruct} The semigroup $\End(\cK)$ of non-degenerate endomorphisms of the $C^*$-algebra $\cK$ is an extension of the group $\Int(\cK)$ of inner automorphisms by the semigroup $\nt$.
\end{thm}
Let $\rho\in\End(\cK)$ be a non-degenerate endomorphism, then $\rho$ acts on the $K$-theory group $K_0(\cK)=\Z$ and respects the order. Thus its action is given by multiplication by an integer $\mmod(\rho)\in \nt$ and this provides a homomorphism of semigroups $\mmod:\End(\cK)\to \nt$ whose kernel is the group $\Int(\cK)$ of inner automorphisms.
It is important to construct explicitly a splitting of this extension and one way to do that is as follows, this also fits with our considerations in characteristic $1$.  We let $\cH:=\ell^2(\Z)$ be the Hilbert space of square integrable sequences of complex numbers. For each integer $k\in \nt$ one considers the isometries $u(k,j):\cH\to \cH$, $0\leq j<k$, which are defined on the canonical orthonormal basis $\delta_n$ of the Hilbert space $\cH:=\ell^2(\Z)$ by
\begin{equation}\label{isoms}
u(k,j)(\delta_n):=\delta_{kn+j} \qqq n\in\Z.
\end{equation}
The maps $\phi(k,j):\Z\to \Z$, $n\mapsto \phi(k,j)(n):= kn+j$ are injective and have disjoint ranges for fixed $k\in \nt$ when $0\leq j<k$. Thus the $u(k,j):\cH\to \cH$ are isometries and their ranges are pairwise orthogonal so that they form a $k$-dimensional  Hilbert space of isometries \ie fulfill the  rules 
\begin{equation}\label{isoms1}
u(k,i)^*u(k,j)=\delta_i^j\qqq i, j \in \{0,\ldots k-1\}, \ \ \sum u(k,i)u(k,i)^*=1.
\end{equation}
Moreover, the composition $\phi(k,j)\circ\phi(k',j')$ is given by $n\mapsto k(k'n+j')+j=kk'n+(kj'+j)$. When written in terms of the isometries $u(k,i)$ one gets the following equality 
\begin{equation}\label{isoms3}
 u(k,j)u(k',j')=u(kk',kj'+j).
\end{equation}
We can then define an action of $\nt$ on the $C^*$-algebra $\cK$ as follows
\begin{lem}\label{isoms2} Let $\cK$ be the $C^*$-algebra of compact operators in the Hilbert space $\cH=\ell^2(\Z)$. The following formula defines an action of the semigroup $\nt$ by endomorphisms of $\cK$
\begin{equation}\label{isoms4}
\fr_k(x):=\sum_{j\in\{0,\ldots k-1\}}u(k,j)\, x\, u(k,j)^*.
\end{equation}
\end{lem}
\proof The multiplicativity of $\fr_k$ follows from \eqref{isoms1} which implies
$$
\fr_k(x)\fr_k(y)=\sum_{i}u(k,i)\, x\, u(k,i)^*\sum_{j}u(k,j)\, y\, u(k,j)^*=\sum_{i}u(k,i)\, xy\, u(k,i)^*=\fr_k(xy).
$$
The equality $\fr_{kk'}=\fr_k \circ\fr_{k'}$ follows from \eqref{isoms3} which implies
$$
\fr_k \circ\fr_{k'}(x)=\sum_{(j,j')}  u(k,j)u(k',j')x \,u(k',j')^*u(k,j)^*=\sum_i u(kk',i)x\,u(kk',i)^*=\fr_{kk'}(x)
$$
since when $(j,j')$ vary in $\{0,\ldots k-1\}\times \{0,\ldots k'-1\}$, one obtains the $kk'$ isometries  $u(kk',i)$ as every element  $i \in \{0,\ldots kk'-1\}$ is uniquely of the form $i=kj'+j$ for $(j,j')\in\{0,\ldots k-1\}\times \{0,\ldots k'-1\}$.\endproof

As a corollary of Lemma \ref{isoms2}, we can think of $\cK$ as an object of the topos $\wnt$, \ie  a sheaf of sets on $\wnt$ and we can use the compatibility of the action of $\nt$ by endomorphisms to view it as a sheaf of involutive algebras and obtain the relevant structure on the stalks of this sheaf. The construction of the stalk given in \eqref{geometricreal} is that of an inductive limit and one needs to handle the distinction between inductive limits taken in the category of $C^*$-algebras with set-theoretic inductive limits. To obtain a $C^*$-algebra from the stalk at a point $\ffp$ of the topos $\wnt$ one needs to take the completion of the set-theoretic stalk. One then has the following corollary of Theorem 5.1 of \cite{Dixmier}:
\begin{prop}\label{propdix} Let $\cK$ be the $C^*$-algebra of compact operators in $\ell^2(\Z)$ endowed with the action of $\nt$ given by Lemma \ref{isoms2}. Then 

$(i)$~$\cK$ is an involutive algebra in the topos $\wnt$.

$(ii)$~The stalk of $\cK$ at a point $\ffp$ of the topos $\wnt$ is a pre-$C^*$-algebra and its completion $\cK_\ffp$ is an infinite separable matroid $C^*$-algebra. 
 
$(iii)$~Let $A$ be an infinite separable matroid $C^*$-algebra, then there exists a point $\ffp$ of the topos $\wnt$, unique up to isomorphism,
 such that $\cK_\ffp$ is isomorphic to $A$.
\end{prop}

Here, by a pre-$C^*$-algebra we mean an involutive algebra (over $\C$) in which the following equality defines a norm whose completion is a $C^*$-algebra
$$
\Vert x\Vert :=\sup \{\vert\lambda \vert\mid \lambda\in \C, \ (x-\lambda)\notin \tilde A^{-1}\}, \  \  \tilde A:= A\oplus \C
$$
where we let $\tilde A:= A\oplus \C$ be the algebra obtained by adjoining a unit. 

The invariant at work in Proposition \ref{propdix} is, as already mentioned above, the $K$-group $K_0(A)$ which is an ordered group called the dimension group of $A$.

%\begin{rem}{\rm Matroid $C^*$-algebras can be non-separable and one may wonder how this fits in the above classification. Let us give a specific example. We let $I$ be a set of arbitrary cardinality and consider the tensor product 
%$$
%\otimes_{\alpha\in I} M_2(\C)_\alpha
%$$
%To remain at the algebraic level we may simply take the inductive limit over finite subsets  $E\subset I$ of the finite tensor products $M(E)=\otimes_{\alpha\in E} M_2(\C)_\alpha.$
%Each of these $C^*$-algebras is isomorphic to $M_n(\C)$ for $n=2^{\vert E\vert}$. For each $\alpha\in I$ one gets a projection $e_\alpha\in M({\alpha})$ given by the diagonal matrix with entries $1$ and $0$ on the diagonal. All these projections are pairwise equivalent $e_\alpha\sim e_\beta$ as one sees using $M({\alpha,\beta})$. Thus even though the size of $M(I):=\varinjlim M(E)$ is  not under control, the dimension group remains the same for any infinite set $I$. 
%}\end{rem}

\end{document}